\documentclass{article}

\usepackage{amsmath}
\usepackage{amssymb} 
\usepackage{latexsym} 
\usepackage{epic,eepic,epsfig,psfrag} 
\usepackage{amscd}  

\newcommand{\QED}{\hspace*{\fill}$\Box$\medskip} 
\newcommand{\eQED}{\vspace{-5.5mm}\hspace*{\fill}$\Box$\medskip}
\def\one{\hbox{1\hskip-2.7pt l}}
\def\smone{{\scriptstyle\rm 1\hskip-2.05pt l}}

\def\rd{{\rm d}}
\def\rT{{\rm T}}

\def\rG{{\rm G}}
\def\p{\phi}
\def\ph{\varphi}
\def\a{\alpha}
\def\b{\beta} 
\def\d{\delta} 

\def\ep{\varepsilon} 
\def\e{\eta} 
\def\g{\gamma}
\def\th{\theta} 
\def\r{\rho}
 
\def\t{\tau} 
 
\def\x{\xi} 
 
\def\n{\nu} 
 
\def\z{\zeta} 
\def\o{\omega} 
\def\l{\lambda}

\def\S{\Sigma}
\def\Si{\Sigma} 
 
\def\X{\Xi} 
\def\Om{\Omega} 
\def\P{\Phi} 
\def\cg{\mathfrak{g}}
\def\cA{{\mathcal A}}
 
\def\cC{{\mathcal C}}

\def\cG{{\mathcal G}}

\def\cL{{\mathcal L}} 
 
\def\cN{{\mathcal N}} 
\def\cm{{\mathfrak m}} 
 
\def\cP{{\mathcal P}} 
\def\cR{{\mathcal R}} 
\def\cS{{\mathcal S}}

\def\cU{{\mathcal U}}

\def\CS{{\mathcal C}{\mathcal S}}
\def\E{{\mathcal E}}
\def\R{{\mathbb R}}

\def\H{{\mathbb H}} 
\def\Z{{\mathbb Z}}
 
\def\half{{\textstyle{\frac 12}}} 
\def\im{{\rm im}\,}
\def\laplace{\Delta} 
\def\Pr{{\bf Proof:}\;} 
\def\st{\: \big| \:}

\def\dt{{\rm d}t}
\def\ds{{\rm d}s}

\def\dph{{\rm d}\p}

\def\dr{{\rm d}r}
\def\pd{\partial}
\def\comp{\circ}

\def\ta{{\tilde\alpha}}

\def\tA{{\tilde A}}

\def\tR{{\tilde R}}

\def\tX{{\tilde\X}}

\def\bep{{\bar\ep}}

\def\la{\langle\,}
\def\ra{\,\rangle}
\def\dvol{\,\rd{\rm vol}}

\def\SU{{\rm SU}}
\def\su{{\mathfrak s\mathfrak u}}

\def\Hom{{\rm Hom}}

\def\tsum{\textstyle\sum}
\def\tint{\textstyle\int}

\def\tr{{\rm tr}}

\newcommand{\winner}[2]{\langle #1{\wedge}#2\rangle}

\newtheorem{dfn}{Definition}[section] 
\newtheorem{lem}[dfn]{Lemma} 
\newtheorem{prp}[dfn]{Proposition} 
\newtheorem{thm}[dfn]{Theorem} 
\newtheorem{rmk}[dfn]{Remark} 
\newtheorem{cor}[dfn]{Corollary}

\begin{document}

\bibliographystyle{plain}

\author{Tomasz S.\ Mrowka, Katrin Wehrheim}

\title{$L^2$-topology and Lagrangians in the space of connections
over a Riemann surface}








\maketitle


\section{Introduction}

Let $\S$ be a closed Riemannian surface and let $\rG$ be a compact Lie group,
whose Lie algebra will be denoted by $\cg$.
We will consider a trivialized $\rG$-bundle over $\S$.
The first goal of this paper is to provide some understanding of the $L^2$-topology
on the space of connections $\cA(\S)=\Om^1(\S;\cg)$ and its quotient by the
gauge group $\cG(\S)=\cC^\infty(\S,\rG)$ acting by bundle isomorphisms.
In particular, we prove a local slice theorem, local connectivity, and uniform local quasiconvexity of the gauge orbits.
The importance of these questions stems from the Yang-Mills path integral over the space of connections, which should naturally be defined using the $L^2$-metric.

The second part of this paper provides some control of the $L^2$-geometry of gauge invariant Lagrangian submanifolds in $\cA(\S)$. These are the natural boundary conditions for a Yang-Mills Floer theory on $3$-manifolds with boundary $\S$ developed in \cite{SW}.
The underlying compactness results for moduli spaces of anti-self-dual instantons with Lagrangian boundary conditions are established in \cite{W bubb} for special Lagrangians arising from handle bodies bounding $\S$. 
Extending the compactness and hence Yang-Mills Floer theory to general gauge invariant Lagrangians requires a weak bound on curvature and local quasiconvexity, which we establish based on a quantitative version of local connectivity of gauge orbits.\\

{\it We thank Stefan Wenger for help with the general analysis of metric spaces, and the meticulous referee for help with the exposition.}

\subsection{Local connectivity and quasiconvexity of gauge orbits}

The action of the gauge group, $u^*A=u^{-1}A u + u^{-1}\rd u$, is a smooth map
$\cG^{1,p}(\S)\times\cA^{0,p}(\S)\to\cA^{0,p}(\S)$ with respect to the
$W^{1,p}$- and $L^p$-topologies for any $p>2$.
The $W^{1,2}$-closure of $\cG(\S)$ however is not a Banach Lie group.
In order to achieve a group structure and a smooth action on the $L^2$-closure of the space of connections $\cA^{0,2}(\S)$, one would have to use the $W^{1,2}\cap L^\infty$-topology on the gauge group.
We will instead fix some $p>2$ and study the gauge orbits in $\cA^{0,p}(\S)$
with respect to the $L^2$-topology.
Hence, in the following we denote by $B_\ep(A_0)\subset \cA^{0,p}(\S)$ the open $L^2$-ball of radius $\ep>0$ around $A_0$.
In Section~\ref{proofs} we prove uniform local quasiconvexity and local pathwise connectedness of the gauge orbits, as stated in Theorem~\ref{thm im Kleinen}, and defined below.

\begin{thm} \label{thm im Kleinen}
For every connection $B\in\cA^{0,p}(\S)$ the gauge orbit
$\cG^{1,p}(\S)^*B$, equipped with the $L^2$-topology, is (i) locally pathwise connected and (ii) uniformly locally quasiconvex.
More precisely:
\begin{enumerate}
\item
$\cG^{1,p}(\S)^*B$ is locally pathwise connected:
Given any $\ep>0$ one can find $\delta>0$
such that for any $A_0,A_1 \in\cG^{1,p}(\S)^*B$  with $\|A_0-A_1\|_{L^2}\leq\delta$ there exists a continuous path
$[0,1]\to B_\ep(A_0)\cap\cG^{1,p}(\S)^*B$, $t\mapsto A_t$ connecting $A_0$ to $A_1$.
\item
The path $t\mapsto A_t$ in (i) can be chosen such that it is smooth as path in $\cA^{0,p}(\S)$, with derivative $\|\partial_t A_t\|_{L^2} \leq C \|A_1-A_0\|_{L^2}$. 
In particular (i) holds with $\delta=\min\{\ep/{2C},\delta_0\}$ for some $\delta_0>0$.
The constants $C, \delta_0$ depend on $[B]\in\cA^{0,p}(\S)/\cG^{1,p}(\S)$.
\item
If $B_0\in\cA^{0,p}(\S)$ is irreducible then there exists an $L^2$-neighbourhood of $B_0$ in $\cA^{0,p}(\S)$, such that on any $L^p$-bounded subset the constant $C$ in (ii) can be chosen uniform.
\end{enumerate}
\end{thm}

The question of a uniform linear relation $\delta=c\ep$ (for $\delta\leq\delta_0$) in Theorem~\ref{thm im Kleinen}~(i)
or a uniform constant $C$ for different gauge orbits in the local convexity (ii) is open for $L^p$-neighbourhoods of reducibles as well as for $L^2$-neighbourhoods.
A positive answer would greatly simplify the proof of local convexity of gauge invariant Lagrangian submanifolds.

The proof of Theorem~\ref{thm im Kleinen} in Section~\ref{proofs} will be based on a subtle $L^2$-local slice theorem explained in Section~\ref{subsec slice}. The remainder of this subsection clarifies the notions of quasiconvexity and local pathwise connectedness, and their relations.

\begin{dfn} \label{dfn path connected}
A topological space $X$ is called {\em locally pathwise connected} if for every open set $\cU\subset X$ and any point $x\in\cU$ the path connected component 
$$
\cP_{\cU,x} :=\bigl\{ y \in \cU \st \exists \gamma\in\cC^0([0,1],\cU) : \gamma(0)=x, \gamma(1)=y \}
$$
is a neighbourhood of $x$.
\end{dfn}


\begin{rmk}
The usual definition of local pathwise connectedness requires a neighbourhood basis of open, pathwise connected sets. This is equivalent to the pathwise connected components $\cP_{\cU,x}$ of all open sets being open; which in turn is equivalent to our Definition~\ref{dfn path connected} above. 
Indeed, note that $\cP_{\cU,x}=\cP_{\cU,y}$ for any $y\in\cP_{\cU,x}$. Hence if $X$ satisfies our definition, then $\cP_{\cU,x}$ is a neighbourhood of $y$ for each $y\in\cP_{\cU,x}$, and hence this pathwise connected component is open.
\end{rmk}

We will deduce pathwise connectivity from the following version of local quasiconvexity with uniform constants. The latter is a general notion for metric spaces, which is closely related to local quasiconvexity as defined by Heinonen \cite[p.57]{Heinonen}.

\begin{dfn} \label{dfn loc q}
Let $(X,d)$ be a metric space.
\begin{enumerate}
 \item 
The {\em length} of a path $\gamma\in\cC^0([0,1],X)$ is 
$$
\ell(\gamma) = \sup\bigl\{\tsum_i d(\gamma(t_i), \gamma(t_{i+1})) \,|\, 0=t_0<t_1 ...< t_k = 1\bigr\} .
$$
\item
$(X,d)$ is {\em locally quasiconvex} if every point $x\in X$ has a neighbourhood $\cN\subset X$ that is quasiconvex. That is, there is a constant $C$ such that any two points $\gamma(0),\gamma(1)\in\cN$ can be joined by a path  $\gamma\in\cC^0([0,1],\cN)$ of length $\ell(\gamma)\leq Cd(\gamma(0),\gamma(1))$. 
\item
$(X,d)$ is {\em uniformly locally quasiconvex} if there exist
constants $\delta>0$ and $C$ such that any two points 
$\gamma(0), \gamma(1)\in X$ with $d(\gamma(0),\gamma(1))\leq\delta$ 
can be joined by a path  $\gamma\in\cC^0([0,1],X)$
of length $\ell(\gamma)\leq Cd(\gamma(0),\gamma(1))$. 
\end{enumerate}
\end{dfn}

\begin{rmk}
In our applications, the paths $\gamma:[0,1]\to \cG^{1,p}(\S)^*B$ in a gauge orbit will be smooth as maps to $\cA^{0,p}(\S)$, and hence 
$\ell(\gamma)=\int_0^1 \| \partial_s \gamma(s) \|_{L^2(\S)}$.
\end{rmk}

In order to justify our notion of uniform local quasiconvexity, we note that the notion of local quasiconvexity above directly implies the following property.
\begin{itemize}
\item[{\em($*$)}]
For all $x\in X$ there exist $\ep>0$ and $C\geq 1$ such that for all $\gamma(0), \gamma(1) \in B_\ep(x)$ there exists a continuous path  $\gamma\in\cC^0([0,1],X)$ of length $\ell(\gamma)\leq Cd(\gamma(0),\gamma(1))$.
\end{itemize}
On the other hand, suppose that $X$ is compact.
Then, firstly, $(*)$ implies local quasiconvexity. (The proof is elementary yet somewhat lengthy.) 
Secondly, $(*)$ is equivalent to uniform local quasiconvexity (as follows from the Lebesgue Lemma).
So for compact metric spaces, local quasiconvexity is equivalent to uniform local quasiconvexity.

In general, neither of the notions of local quasiconvexity and uniform local quasiconvexity in Definition~\ref{dfn loc q} implies the other. However, either of them (as well as $(*)$ above) implies local connectivity. We only prove the part that is relevant in our setting.

\begin{prp}
If a metric space is uniformly locally quasiconvex, then it is locally pathwise connected as topological space.
\end{prp}
\Pr
Let an open set $\cU\subset X$ and a point $x\in\cU$ be given.
Since $\cU$ is a neighbourhood of $x$, it contains a metric ball $B_\ep(x)$ for some $\ep>0$.
Now choose $r>0$ such that $r\leq\delta$ and $Cr < \ep$ with the constants from the uniform local
quasiconvexity. 
Now any point $y\in B_r(x)$ can be connected to $x$ by a path $\gamma$ of length $\ell(\gamma)\leq Cr<\ep$.
By the definition of the length, this path must be entirely contained in $B_\ep(x)\subset\cU$. (Indeed, $\gamma(t')\in X\setminus B_\ep(x)$ would imply
$\ell(\gamma)\geq d(x,\gamma(t')) + d(\gamma(t'),y) > \ep + 0$.)
Hence the path connected component of $\cU$ containing $x$ also contains the ball $B_r(x)$, and hence is a neighbourhood of $x$.
\QED

\subsection{The $\mathbf{L^2}$-local slice theorem} \label{subsec slice}

The study of the $L^2$-topology on the moduli space of connections $\cA^{0,p}(\S)/\cG^{1,p}(\S)$ hinges on a local slice theorem which provides generalized orbifold charts.
Here we fix $p>2$ and work with the space of connections
$\cA^{0,p}(\S)=L^p(\S;T^*\S\otimes\cg)$ and the gauge group $\cG^{1,p}(\S)=W^{1,p}(\S,\rG)$. 
A slice of the gauge action at a connection $A_0\in\cA^{0,p}(\S)$ is the $L^2$-orthogonal complement to the gauge orbit,
\begin{align*}
S_{A_0}&:= 
\bigl\{ A=A_0+a\in\cA^{0,p}(\S) \st \rd_{A_0}^*a = 0 \bigr\}.
\end{align*}
Here $\rd_{A_0}^* : L^p(\S,\rT^*\S\otimes\cg) \to W^{-1,p}(\S,\cg):=W^{1,p'}(\S,\cg)^*$ 
with ${p^{-1}+p'^{-1}=1}$
is the dual operator of $\rd_{A_0}: W^{1,p'}(\S,\cg)\to  L^{p'}(\S,\rT^*\S\otimes\cg)$.
The content of the local slice Theorem~\ref{thm local slice} below is that the gauge orbits 
through $\cS_{A_0}$ cover an $L^2$-neighbourhood of $A_0$.
This provides generalized orbifold chart for $L^2$-neighbourhoods on 
$\cA^{0,p}(\S)/\cG^{1,p}(\S)$. The symmetry group will be the stabilizer
$$
{\rm Stab}(A_0):=\bigl\{ g\in\cG^{1,p}(\S) \st g^*A_0 = A_0 \bigr\},
$$
which is a nondiscrete but compact Lie group for reducible $A_0$.
We will denote the $L^2$-balls in the local slice of radius $\ep>0$ by
\begin{align*}
S_{A_0}(\ep)&:= 
\bigl\{ A=A_0+a\in\cA^{0,p}(\S) \st \rd_{A_0}^*a = 0 , \|a\|_{L^2}<\ep \bigr\} .
\end{align*}
These are invariant under the
action $A_0+a \mapsto g^*(A_0+a)=A_0+g^{-1}a g$ of the stabilizer
$g\in{\rm Stab}(A_0)$.

\begin{thm} \label{thm local slice} {\bf (Local Slice Theorem)}
For every $A_0\in\cA^{0,p}(\S)$ there are constants ${\ep,\d>0}$ 
such that the map
\begin{equation}\label{map}
\cm:\; \begin{aligned}
\bigl(\cS_{A_0}(\ep) \times \cG^{1,p}(\S)\bigr) / {\rm Stab}(A_0) 
\; & \to & \cA^{0,p}(\S) \\
[(A_0+a,u)] \;& \mapsto & u^{-1\;*}(A_0 + a)
\end{aligned}
\end{equation}
is a diffeomorphism onto its image, which contains an $L^2$-ball,
$$
B_\d(A_0):=
\bigl\{ A\in\cA^{0,p}(\S) \st \|A-A_0\|_{L^2}<\d \bigr\} 
\subset \im\cm .
$$
\end{thm}

\begin{cor}  \label{cor local slice}
Any $A\in\cA^{0,p}(\S)$ with $\|A-A_0\|_{L^2}<\d$ is gauge equivalent
to a connection in the local slice through $A_0$, that is
$\rd_{A_0}^*( u^*A - A_0 )=0$ for some gauge transformation $u\in\cG^{1,p}(\S)$.
Moreover, the gauge transformation $u$ is unique up to ${\rm Stab}(A_0)$.
\end{cor}

The proof of the local slice theorem~\ref{thm local slice} for base manifolds of dimension $n\geq 3$ can be found in \cite{M}. We restrict our proof in Section~\ref{proofs} to the case $n=2$, which is somewhat more complicated since the simple estimate $\|f g\|_{W^{-1,n}} \leq C \|f\|_{L^n}\|g\|_{L^n}$ has to be replaced by the div-curl Lemma from harmonic analysis (see e.g.\ \cite{Taylor}).
In Section~\ref{est} we prove a weaker version, Lemma~\ref{lem div curl}, and extend it to our gauge theoretic settings.

\subsection{Yang-Mills Floer theory on 3-manifolds with boundary}

The natural symplectic form on the space of connections $\cA(\S)=\Om^1(\S;\cg)$ is
\begin{equation} \label{omega}
\o(\a,\b) := \int_\Si \la \a \wedge \b \ra 
\qquad\quad\text{for}\;\; \a,\b\in\Om^1(\S;\cg) .
\end{equation}
Here the values of the differential forms are paired by the inner product
$\la \cdot , \cdot \ra$ on $\cg$. 
This symplectic form appears in the work of Atiyah and Bott \cite{AB}, who observed that the moduli space of flat connections on $\S$ is on the one hand homeomorphic to the compact representation space $\cR_\S = \Hom(\pi_1(\S),\rG)/\rG$, and can on the other hand be viewed as the symplectic quotient of the gauge action on the infinite dimensional space of connections,
$$
\cR_\S \;\cong\;  \{ A\in \cA(\S) \st F_A=0\} /\cG(\S) 
\;=\; \cA(\S) /\hspace{-1mm}/ \cG(\S) .
$$
This is because the moment map with respect to the symplectic structure \eqref{omega} is given by the curvature,
$\cA(\S) \rightarrow \Om^0(\S;\cg)=\rT_1\cG(\S)$, $A \mapsto *F_A = * (\rd A + \tfrac 12 [A \wedge A] )$.

The symplectic form \eqref{omega} also naturally appears in Yang-Mills field theory:
The 
Yang-Mills functional on the space of connections on a $4$-manifold $X$ is
$\tA\mapsto \half \int_X \la F_\tA \wedge F_\tA \ra$. For a compact $4$-manifold with boundary $\pd X=Y$ it equals to 
$$
\half \int_X \la F_\tA \wedge F_\tA \ra =
\half \int_Y \la \tA \wedge \bigl( F_\tA - \tfrac 16 [\tA\wedge \tA] \bigr) \ra  =: \CS(\tA|_Y) 
\qquad\text{for}\; \tA\in\cA(X),
$$
which defines the Chern-Simons functional $\CS:\cA(Y)\to\R$ for a compact $3$-manifold $Y$. 
It only descends to a multivalued functional on the moduli space of connections by the gauge action, but its differential defines the gauge invariant Chern-Simons $1$-form
$$
\l_A (\a) :=  \int_Y \la F_A \wedge \a \ra 
\qquad\quad\text{for}\;\; \alpha\in\rT_A\cA(Y)=\Om^1(Y;\cg).
$$
Instanton Floer theory for a closed $3$-manifold $Y$ is the Morse theory for this closed $1$-form on $\cA(Y)/\cG(Y)$, as developed by Floer \cite{F1}.
The Chern-Simons $1$-form is also well defined and gauge invariant on a $3$-manifold with boundary $\pd Y=\S$, but it is no longer closed. 

In fact, its differential is the symplectic form (\ref{omega}):
For $\a,\b\in\Om^1(Y;\cg)=\rT_A\cA(Y)$
\begin{align*}
\rd\l(\a,\b) 
&\;=\; \int_Y \la \rd_A\a\wedge\b \ra - \int_Y \la \rd_A\b\wedge\a \ra
\;=\; \int_{\S} \la \a|_\S \wedge \b|_\S \ra .
\end{align*}
To render $\l$ closed, it is natural to pick a Lagrangian submanifold $\cL\subset\cA(\S)$ and restrict the $1$-form to the space of connections $\cA(Y,\cL):=\{A\in\cA(Y)\,|\,A|_\S\in\cL\}$ with Lagrangian boundary value.
(If $\cL\subset\cA(\Si)$ is any submanifold, then the closedness of $\l|_{\cA(Y,\cL)}$
is equivalent to $\o|_\cL\equiv 0$, and the maximal such submanifolds are precisely
the Lagrangian submanifolds.)
If we in addition assume $\cL$ to be gauge invariant then $\l$ descends to a closed $1$-form on $\cA(Y,\cL)/\cG(Y)$, from which instanton Floer theory for the pair $(Y,\cL)$ is developed in \cite{SW}.
Its critical points are the flat connections with Lagrangian boundary value $\{ F_A=0, A|_\S\in\cL\}/\cG(Y)$, and trajectories defining the differential are gauge equivalence classes of anti-self-dual connections on $\R\times Y$ with Lagrangian boundary conditions, i.e.\ 
connections $\X\in\cA(\R\times Y)$ that solve the boundary value problem
\begin{equation}\label{bvpi}
F_\X + * F_\X = 0, \qquad\qquad
\X|_{\{s\}\times\Si} \in\cL \quad\forall s\in\R .
\end{equation}
In local coordinates near the boundary, with a metric $\ds^2+\dt^2+g_\S$ on ${\R\times[0,\delta)\times\S}$, the system \eqref{bvpi} for the connection $\Xi=\P\ds+\Psi\dt+A$ can be rewritten in terms of maps $A:\R\times[0,\delta)\to\cA(\S)$ and 
$\P,\Psi:\R\times[0,\delta)\to\cC^\infty(\S,\su(2))\cong\rT_\smone\cG(\S)$ that satisfy
\begin{equation} \label{eq 1} 
\left\{ \begin{aligned}
( \pd_s A - \rd_A\P ) + * ( \pd_t A - \rd_A\Psi ) &= 0 ,\\
\pd_s \Psi - \pd_t\P + [\P,\Psi] + * F_A  &=0 , 
\end{aligned} \right.
\qquad \qquad
A(s,0) \in \cL \quad \forall s\in\R .
\end{equation}
These are the symplectic vortex equations \cite{CGMS} for the $\cG(\S)$-action on $\cA(\S)$. Indeed, $\P\mapsto\rd_A\P$ is the infinitesimal action,
$A\mapsto *F_A$ is the moment map, and the Hodge operator $*$ of any metric $g_\S$ is a complex structure on $\cA(\S)$ compatible with the symplectic form \eqref{omega}.
However, unlike the finite dimensional setting in \cite{CGMS}, the compactness of moduli spaces of solutions of \eqref{eq 1} is far from clear since the boundary condition is a combination of a first order condition (for a first order equation!) and nonlocal conditions.
This is due to the fact that the gauge invariant Lagrangians lie in the subset of flat connections\footnote{
This follows directly if $[\cg,\cg]=\cg$ as for any Lie groups with discrete center. For general Lie groups it is a natural restriction.
} and hence are determined by a (singular) Lagrangian $\cL/\cG(\S)\subset\cR_\S$ in the representation space.
So the boundary condition $A(s,0) \in \cL$ can equivalently be rewritten as $F_{A(s,0)}=0$ and a finite number ($\half\dim\cR_\S$) of holonomy conditions for $A(s,0)$.
The key to the proof of elliptic estimates for \eqref{bvpi} in \cite{W elliptic} is to view $\pd_s A  + * \pd_t A = *\rd_A\Psi  - \rd_A\P$ as perturbed holomorphic curve $A:\R\times[0,\delta)\to\cA^{0,p}(\S)$ in a symplectic Banach space.
In fact, the bubbling analysis in \cite{W bubb} indicates that solutions of \eqref{eq 1} truly behave like holomorphic curves near the boundary, in that they concentrate energy and develop singularities along slices $\{s\}\times \S \subset \R\times Y$ rather than at points, as anti-self-dual instantons on closed $4$-manifolds do.
A complete geometric description of the bubbles arising from rescaling near the singularity
(as holomorphic disks in $\cA(\S)$) is still open but \cite{W bubb} succeeded in proving compactness of moduli spaces for $G=\SU(2)$ and the special Lagrangian submanifolds 
$$
\cL_H:=\{ \tA|_\S \,|\, \tA\in\cA(H), F_\tA = 0 \} \subset \cA(\S)
$$ 
that arise from a handle body $H$ bounding $\S=\pd H$. The two crucial ingredients are an energy quantization for the bubbles and removal of singularity for solutions of \eqref{eq 1} on the complement of a slice $\{(s,0)\}\times \S$. Again, the proofs use symplectic techniques for holomorphic curves with boundary on $\cL\subset\cA^{0,p}(\S)$. They hinge on certain bounded geometry of the Lagrangian, which are evident for compact Lagrangians in finite dimensional symplectic manifolds, far from obvious in the infinite dimensional Banach space setting, and in the special case of $\cL_H$ can be replaced by a subtle extension theorem from $W^{1,2}(\S,\SU(2))$ to $W^{1,3}(H,\SU(2))$, due to Hardt-Lin \cite{HrL} in this borderline Sobolev case.

In Section~\ref{sec Lag} we provide the required geometric control for a general class of gauge invariant  Lagrangian submanifolds $\cL\subset \cA(\S)$ (see Definition~\ref{dfn L} for details).
In particular, we show uniform local quasiconvexity, which is then used in Section~\ref{sec CS} to define a local Chern-Simons functional for paths in $\cA(\S)$ with ends on $\cL$.
Section~\ref{sec bubb} builds on these results to prove the compactness of moduli spaces of solutions of \eqref{eq 1} for general gauge invariant Lagrangians.
We explain the overall philosophy in the next section.

\subsection{Curvature and local quasiconvexity of gauge invariant Lagrangians}

The first main step in the proof of compactness for \eqref{eq 1} is an energy quantization: {\em If in a sequence of solutions the energy density blows up at a point, then it also concentrates a minimal quantum of energy at that point.} For holomorphic curves with Lagrangian boundary values $u:(\Omega,\partial\Omega)\to (M,L)$, this can be proven without identifying the bubbles as spheres or disks, but just using a mean value inequality for the energy density $e=|\rd u|^2$, which arises from a nonlinear subharmonicity (see e.g.\ \cite{W mean} -- here stated with the positive definite Laplace operator)
$$
\laplace e \leq  A + a e^{2} , \qquad
\tfrac\pd{\pd\n}\bigr|_{\pd\Omega} e \leq B + b e^{3/2} .
$$
For anti-self-dual connections, the $4$-dimensional energy density $\half |F_\X|^2=|\pd_s A|^2 + |F_A|^2$ satisfies the corresponding nonlinear bound on the Laplacian, but due to the nonlocal boundary conditions, one has no control over the normal derivative. 
Viewing $A:\R\times[0,\delta)\to\cA(\S)$ as holomorphic curve, one is led to trying the $2$-dimensional energy density $e=\|\pd_s A\|_{L^2(\Sigma)}^2 +\|F_{A}\|_{L^2(\Sigma)}^2$. Its Laplacian causes some difficulties (that are dealt with in \cite{W bubb}), but the normal derivative simply behaves just like a holomorphic curve, 
$$
- \half \tfrac\pd{\pd t} e \bigr|_{t=0}
\;\leq\;  C \| \pd_s A \|_{L^2(\S)}^2 
     + \omega( \pd_s A , \pd_s^2 A ) 
\;\leq\;  B + b \|\pd_s A \|_{L^2(\S)}^{3/2}    .
$$
Indeed, if $\cL\subset \cA(\S)$ was a Lagrangian subspace then $\omega( \pd_s A , \pd_s^2 A )=0$ since $\partial_s A,\partial_s^2 A \in \rT\cL$. If $\cL$ was a compact Lagrangian, then $\omega( \pd_s A , \pd_s^2 A )\leq C |\pd_s A|^3$ with a constant C arising from the curvature of $\cL$. In our nonlinear infinite dimensional setting, it is not clear whether curvature is uniformly bounded, but the essential estimate can be shown directly.

\begin{lem} \label{lem lin ext i}
There is a constant $C_{\rT\cL}$
such that any smooth path $A:\R\to\cL$ satisfies
\begin{align*}
\int_\S \la \pd_s A(0) \wedge \pd_s\pd_s A(0) \ra 
&\leq C_{\rT\cL} \bigl\| \pd_s A(0) \bigr\|_{L^2(\Si)}^3  .
\end{align*}
\end{lem}

For a Lagrangian $\cL_H$ arising from a handle body $\partial H=\S$, this was proven in \cite{W bubb} using the extension to flat connections $\tA\in\cA_{\rm flat}(H)$, which satisfy $\rd_\tA\pd_s\tA=0$ and hence
$$
 \int_\Si \la \pd_s A \wedge \pd_s^2 A \ra  
= \int_H \rd \la \pd_s \tA \wedge \pd_s^2 \tA \ra  
=  \int_H \la \pd_s \tA \wedge [\pd_s\tA \wedge \pd_s\tA] \ra  
\leq \|\pd_s\tA\|_{L^3(H)}^3 .
$$
So Lemma~\ref{lem lin ext i} follows from constructing extensions with
$\|\pd_s\tA\|_{L^3(H)} \leq C_H \|\pd_s A\|_{L^2(\S)}$.
For a general gauge invariant Lagrangian we 
restate and prove this result as Lemma~\ref{lem lin ext},
using the local slice Theorem~\ref{thm local slice} and a div-curl Lemma from Section~\ref{est}.\\

The other main step in the proof of compactness for \eqref{eq 1} is a removal of singularity for solutions on the complement of a slice $\{(s,0)\}\times \S$ with finite energy. Again, the proof in \cite{W bubb} for the special case $\cL_H$ proceeds by combining gauge theoretic analysis with holomorphic curve techniques. 
In particular, one needs an analogue of the local symplectic action: 

If $\Xi\in\cA(B^*\times\S)$ is an anti-self-dual connection over the punctured ball $B^*=B_r \setminus\{0\}\subset \R^2$, and in polar coordinates there is a gauge $\X=A+R\dr$ with no $\dph$-component (i.e.\ there is no holonomy around the singularity), then one can express the energy on annuli as difference of integrals over the boundary, just as for a punctured holomorphic curve $u:B^*\to (M,\omega)$. Indeed, we have the same formula for energies
$$
\half \int_{(B_\rho\setminus B_\delta)\times\S} |F_\X|^2  = \cC(A(\delta,\cdot)) - \cC(A(\rho,\cdot)) ;
\qquad
\int_{B_\rho\setminus B_\delta} u^*\omega  = \cA(u(\delta,\cdot)) - \cA(u(\rho,\cdot)) ,
$$
where the boundary integrals depend on a choice of gauge resp.\  local primitive $\omega=\rd\lambda$,
$$
\cC(A:S^1\to\cA(\S)) = - \half\int_0^{2\pi} \int_\S \la A(\p) \wedge \pd_\p A(\p) \ra \;\dph ;
\qquad
\cA(u:S^1\to M) = - \int_0^{2\pi} u^*\lambda .
$$
The local symplectic action is a canoncial definition of $\cA(u:S^1\to M)$ for all sufficiently short loops, given by a unique local choice of primitive.
There is a similar construction for holomorphic curves with Lagrangian boundary values, where the local symplectic action of a path with endpoints on the Lagrangian is given by connecting the endpoints within the Lagrangian to define a loop.
The crucial analytic ingredient to the removal of holomorphic singularities then is the isoperimetric inequality for sufficiently short loops,
$$
|\cA(u: S^1\to M) | \;\leq\; C 
\Bigl( \int_0^{2\pi} | \pd_\p u | \,\dph \Bigr)^2 .
$$
For anti-self-dual connections $\Xi=A + R\dr + \Phi\d\phi \in\cA(D^*\times\S)$ over the punctured half ball $D^*=B_r \setminus\{0\}\subset \H^2$ with Lagrangian boundary condition $A|_{\pd\H^2}\in\cL$, there always exists a gauge $\X=A+R\dr$ without angular component. As for holomorphic curves, finite energy implies that the paths 
$A(r,\cdot):[0,\pi]\to\cA(\S)$ become short for $r\to 0$. However, this only holds in the $L^2$-norm on $\cA(\S)$ that enters in the energy.
Due to the subtlety of the gauge action in the $L^2$-topology, we cannot assume that the $L^2$-closure $\cL\subset\cA^{0,2}(\S)$ is a Hilbert submanifold; instead we are throughout working with Banach submanifolds $\cL\subset\cA^{0,p}(\S)$ for $p>2$. 
Nevertheless, in order to define a local Chern-Simons functional we now need to know that $L^2$-close points on the $L^p$-submanifold $\cL$ can be connected within the Lagrangian. Moreover, in order to obtain an isoperimetric inequality, the length of the connecting path must be linearly bounded by the distance of the endpoints. In other words, we precisely need the following local convexity with uniform constants.

\begin{lem}  \label{lem ext i}
There are constants $C_\cL$ and $\d_\cL>0$ such that for all 
$A(0), A(1)\in\cL$ with $\|A(0)-A(1)\|_{L^2}\leq\d_\cL$ there exists a smooth connecting path 
$A:[0,1]\to \cL$ such that
$$
 \|\pd_s A (s) \|_{L^2} \leq C_\cL \|A(0)-A(1)\|_{L^2} 
\qquad \forall s\in[0,1].
$$
\end{lem}

We
restate and prove this result in Lemma~\ref{lem ext} based on the local convexity of gauge orbits and assuming that $\cL/\cG(\S)$ is compact. 
We moreover compare our construction with the previously defined local Chern-Simons functional \cite{W bubb} for the special Lagrangians $\cL_H$ arising from handle bodies, employing the following nontrivial extension property:
For all $A(0), A(1)\in\cL_H$ there exist
$\tA(0),\tA(1)\in\cA_{\rm flat}(H)$ with $A(i)=\tA(i)|_\S$ 
such that with a uniform constant $C_H$
$$
\|\tA(0)-\tA(1)\|_{L^3(H)} \leq C_H \|A(0)-A(1)\|_{L^2(\S)} .
$$

\section{Div-Curl Lemmas}
\label{est}

The aim of this section is to establish some estimates that will be crucial for the proofs in the subsequent section. The first is the so-called div-curl lemma from harmonic analysis 
(see e.g.\ \cite{Taylor}).
Here we give an alternative proof of a weaker version (using $W^{1,2}(\Sigma)^*$
rather than BMO spaces) and extend it to our gauge theoretic settings.
Here, as in the introduction, we fix $\S$ to be a closed Riemannian surface, 
let $\rG$ be a compact Lie group with Lie algebra $\cg$, and consider a trivialized $\rG$-bundle over $\S$.
The div-curl lemmas will later be applied to both gauge transformations $h:\S\to\rG$ and 
infinitesimal gauge transformations $h:\S\to\cg$. 
For that purpose we can identify $\rG$ with a subgroup of $\R^{N\times N}$ such that the
Lie bracket is replaced by the commutator $[\x,h]=\x h - h \x$ 
and the inner product becomes the trace $\la\x,\e\ra = -\tr(\x\e)$.
%
%
We will be working with the space of smooth connections $\cA(\S)=\cC^\infty(\S;T^*\S\otimes\cg)$ as well as with the space of $L^p$-connections $\cA^{0,p}(\S)=L^p(\S;T^*\S\otimes\cg)$ for a fixed $p>2$.
Moreover, for any $1<q<\infty$ we denote by $\|\cdot\|_{W^{-1,q}}$ 
the norm on the dual space of $W^{1,q'}(\S)$,
where $q'$ is the dual exponent given by ${q^{-1}+q'^{-1}=1}$.

\begin{lem}{\bf ( Div-Curl Lemma ) } \label{lem div curl}
For every smooth connection $A_0\in\cA(\S)$ there is a constant $C$
such that for all $f,g,h\in\cC^\infty(\S,\R^{N\times N})$
$$
\left| \int_\S \tr\bigl( f \cdot  \rd_{A_0} g \wedge \rd_{A_0} h \bigr) \right|
\;\leq\; 
C \bigl( \|f\|_{L^2} + \| \rd_{A_0} f \|_{L^2} \bigr) 
\| \rd_{A_0} g \|_{L^2} \| \rd_{A_0} h \|_{L^2} .
$$
%
%
%
%
\end{lem}
\Pr
We pick a compact Riemannian $3$-manifold $Y$ that bounds ${\pd Y=\S}$ and 
extend the connection $A_0$ to $\tA_0\in\cA(Y)$.
Next, we fix a complement $\tilde K\subset W^{1,3}(Y,\R^{N\times N})$ of
$\{\tilde f\in \ker\rd_{\tA_0} |\; \tilde f|_{\pd Y}=0 \}$.
Then the solutions $\tilde f\in\tilde K$ of the Dirichlet problem 
${\rd_{\tA_0}^*\rd_{\tA_0}\tilde f=0}$, $\tilde f|_{\pd Y}=f$
are unique and define a bounded linear extension map 
$W^{1,2}(\S,\R^{N\times N})\to W^{1,3}(Y,\R^{N\times N})$, $f\mapsto\tilde f$.
(See e.g.\ \cite[Lemma~3.2]{W bubb} for the continuity.)
The extensions satisfy 
$$
\| \tilde f \|_{W^{1,3}(Y)}\;\leq\; C_1 \| f \|_{W^{1,2}(\S)} 
\;\leq\; C_2 \bigl(  \|f\|_{L^2(\S)} + \| \rd_{A_0} f \|_{L^2(\S)} \bigr).
$$
Without loss of generality we can assume that $g$ and $h$ lie in 
a complement $K\subset W^{1,2}(\S,\R^{N\times N})$ of $\ker\rd_{A_0}$.
Then their extensions satisfy 
$$
\| \tilde g \|_{W^{1,3}(Y)}\;\leq\; C_3 \| \rd_{A_0} g \|_{L^2(\S)} , \qquad
\| \tilde h \|_{W^{1,3}(Y)}\;\leq\; C_3 \| \rd_{A_0} h \|_{L^2(\S)} .
$$
Using these extensions we can estimate for any $f\in\cC^\infty(\S,\R^{N\times N})$
and $g,h\in K$
\begin{align*}
& \biggl| \int_\S \tr\bigl( f \cdot \rd_{A_0} g \wedge \rd_{A_0} h \bigr) \biggr| \\
&= \biggl| 
\int_Y \tr\bigl( \rd_{\tA_0} \tilde f 
           \wedge  \rd_{\tA_0} \tilde g \wedge \rd_{\tA_0} \tilde h \bigr) \\
&\qquad + \int_Y \tr\bigl( \tilde f \cdot
\bigl( [F_{\tA_0}, \tilde g] \wedge \rd_{\tA_0} \tilde h 
-  \rd_{\tA_0} \tilde g \wedge [ F_{\tA_0} , \tilde h ]  \bigr) \bigr)
\biggr| \\
&\leq C_4 \| \tilde f \|_{W^{1,3}(Y)} \, \| \tilde g \|_{W^{1,3}(Y)} \, \| \tilde h \|_{W^{1,3}(Y)} \\
&\leq C_4 C_2 C_3^2 \bigl( \|f\|_{L^2(\S)} + \| \rd_{A_0} f \|_{L^2(\S)} \bigr)
 \| \rd_{A_0} g \|_{L^2(\S)} \, \| \rd_{A_0} h \|_{L^2(\S)} .
\end{align*}
Here all constants $C_1$ to $C_4$ only depend on $A_0$ and the choice of 
its extension $\tA_0$ (which also determines the curvature $F_{\tA_0}$).
%
%
%
%
\QED

The remainder of this section draws more gauge-theoretic conclusions from the above div-curl lemma.
We first establish a simple estimate that we will need repeatedly.

\begin{lem}\label{lem:c0}
For every $A_0\in\cA^{0,p}(\S)$ with $p>2$ and $q\geq 2$ there is a constant $C_0$ such that
\begin{equation} \label{c0}
\|\rd_{A_0}\x\|_{L^q} \leq C_0 \|\rd_{A_0}^* ( \rd_{A_0}\x) \|_{W^{-1,q}} \qquad\forall \xi\in W^{1,q}(\Sigma,\cg)
\end{equation}
\end{lem}
\Pr
The operator $\rd_{A_0}^*: L^{q}(\Sigma,\rT^*\Sigma\otimes\cg) \to W^{-1,q}(\Sigma,\cg)$
is well defined since the Sobolev embedding $W^{1,q'}(\S)\hookrightarrow L^r(\S)$ with
$$
\tfrac 1{r} = 1 - \tfrac 1p - \tfrac 1q  > \tfrac 12 - \tfrac 1q \geq 0
$$ 
ensures the Sobolev multiplication 
$L^q\cdot L^p \cdot L^r \subset L^1$.
The restricted operator $\rd_{A_0}^*|_{\im\rd_{A_0}}$ is injective 
since every element of the kernel has the form $\rd_{A_0}\xi$ 
with $\xi\in W^{1,q}(\Sigma,\cg)\subset W^{1,q'}(\Sigma,\cg)$ 
and satisfies
$\int_\S \la\rd_{A_0}\xi , \rd_{A_0}\psi\ra = 0$ for all $\psi\in W^{1,q'}$,
which with $\psi=\xi$ implies $\rd_{A_0}\xi = 0$.
The cokernel of $\rd_{A_0}^*|_{\im\rd_{A_0}}$ equals to the cokernel of the Laplacian $\rd_{A_0}^*\rd_{A_0}$, which is finite dimensional.
Hence $\rd_{A_0}^*|_{\im\rd_{A_0}}$ is a Fredholm operator and \eqref{c0} follows from the injectivity.
\QED

\begin{cor} \label{cor div curl} \hspace{2mm} \\
\vspace{-5mm}
\begin{enumerate}
\item
For every $p>2$ there is a constant $C$ such that for 
every weakly flat\footnote{
We say $A\in\cA^{0,p}(\S)$ is weakly flat if 
$\int_\S \la \rd_{A_0}\xi \wedge \rd_{A_0}\eta \ra =0$ 
for all $\xi,\eta\in\cC^\infty(\S,\cg)$.
}
connection $A_0\in\cA^{0,p}_{\rm flat}(\S)$
and any $\x\in \Om^0(\S;\cg)$ and $\a_1,\a_2\in\Om^1(\S;\cg)$
\begin{align*}
\left| \tint_\S \la \x \,,\, [ \a_1 \wedge \a_2 ] \ra \right|
&\leq
C \bigl( \|\x\|_{L^2} + \| \rd_{A_0} \x \|_{L^2} \bigr)
\prod_{i=1,2}\bigl( \| \a_i \|_{L^2} + \| \rd_{A_0}\a_i \|_{L^{\frac{2p}{p+2}}} \bigr) , \\
\left| \tint_\S \la \x \,,\, [ \a_1 \wedge \a_2 ] \ra \right|
&\leq
C \bigl( \|\x\|_{L^2} + \| \rd_{A_0} \x \|_{L^2} \bigr)
\prod_{i=1,2}\bigl( \| \a_i \|_{L^2} + \| \rd_{A_0}^*\a_i \|_{L^{\frac{2p}{p+2}}} \bigr) .
\end{align*}
\item
For every $p>2$ and connection $A\in\cA^{0,p}(\S)$
there is a constant $C$ (depending continuously on $A$) such that for 
any $\a_1,\a_2\in\Om^1(\S;\cg)$
\begin{align*}
\| [ \a_1 \wedge \a_2 ] \|_{W^{-1,2}}
&\leq
C \prod_{i=1,2}\bigl( \| \a_i \|_{L^2} + \| \rd_{A}\a_i \|_{W^{-1,p}} \bigr) ,\\
\| [ \a_1 \wedge \a_2 ] \|_{W^{-1,2}}
&\leq
C \prod_{i=1,2}\bigl( \| \a_i \|_{L^2} + \| \rd_{A}^*\a_i \|_{W^{-1,p}} \bigr) .
\end{align*}
\end{enumerate}
\end{cor}

\noindent
{\bf Proof of Corollary~\ref{cor div curl} : }
First note that it suffices to prove the first estimate in both (i) and (ii).
The second estimates then follows by applying the first to
$[\a_1\wedge\a_2]=[*\a_1\wedge *\a_2]$
and noting that $\|\rd_{A_0}*\a_i\|=\|\rd_{A_0}^*\a_i\|$.

In (i) the action 
$(A_0,\x,\a_1,\a_2) \mapsto (u^*A_0, u^{-1}\x u, u^{-1}\a_1 u, u^{-1}\a_2 u)$
of the gauge group $\cG(\S)\ni u$ preserves the inequality,
and we will see that the optimal constant varies continuously 
with $A_0$ in the $L^p$-norm.
It thus suffices to prove the estimate on an $L^p$-neighbourhood 
of any given smooth flat connection $A_{0}$. Then finitely many of these 
cover the compact quotient
$\cA^{0,p}_{\rm flat}(\S)/\cG^{1,p}(\S)\cong \cA_{\rm flat}(\S)/\cG(\S)$.
So we fix $A_0\in\cA^{0,p}_{\rm flat}(\S)$ 
and in the following allow the constants to depend on $A_0$.

We use the Hodge decomposition with respect to the flat connection $A_0$ to
write $\a_i=\rd_{A_0}\z_i + *\rd_{A_0}\e_i + \b_i$ with 
harmonic $\b_i\in\ker\rd_{A_0}\cap\ker\rd_{A_0}^*$.
The harmonic $1$-forms form a finite dimensional space, so for any $p>2$ 
we have a constant $C_1$ such that 
$\|\b_i\|_{L^p}\leq C_1\|\b_i\|_{L^2}\leq C_1 \|\a_i\|_{L^2}$.
Moreover, we have as in Lemma~\ref{lem:c0}
$$
\|*\rd_{A_0}\e_i\|_{L^p} 
\leq C_0 \|\rd_{A_0}^*\rd_{A_0}\e_i \|_{W^{-1,p}}
= C_0 \|\rd_{A_0}\a_i\|_{W^{-1,p}} .
$$
Finally, we also have $\|\rd_{A_0}\z_i\|_{L^2}\leq \|\a_i\|_{L^2}$.
Now we use the H\"older inequality for $L^p\cdot L^2 \hookrightarrow L^{\frac{2p}{p+2}}$,
the Sobolev embedding $W^{1,2}(\S) \hookrightarrow L^{\frac{2p}{p-2}}(\S)$,
and the constant $C$ from Lemma~\ref{lem div curl} to estimate 
\begin{align*}
&\left| \int_\S \la \x \,,\, [ \a_1 \wedge \a_2 ] \ra \right| \\
&\leq
\left| \int_\S \la \x \,,\, [ ( \b_1 + *\rd_{A_0}\e_1 ) \wedge \a_2 ] \ra \right|
+ \left| \int_\S \la \x \,,\, [ \rd_{A_0}\z_1 \wedge ( \b_2 + *\rd_{A_0}\e_2) ] \ra \right| \\
&\quad
+ \left| \int_\S \la \x \,,\, 
\bigl( \rd_{A_0}\z_1 \wedge \rd_{A_0}\z_2 + \rd_{A_0}\z_2 \wedge \rd_{A_0}\z_1 \bigr)
\ra \right| \\
&\leq 
\|\x\|_{L^{\frac{2p}{p-2}}} \bigl( \| \b_1 + *\rd_{A_0}\e_1 \|_{L^p} \|\a_2\|_{L^2}
+ \|\rd_{A_0}\z_1\|_{L^2} \|\b_2 + *\rd_{A_0}\e_2) \|_{L^p} \bigr) \\
&\quad + 2 C \bigl(\|\x\|_{L^2} + \|\rd_{A_0}\x\|_{L^2} \bigr) 
  \|\rd_{A_0}\z_1\|_{L^2} \|\rd_{A_0}\z_2 \|_{L^2} \\
&\leq C_3 \|\x\|_{W^{1,2}}
\bigl( \| \a_1 \|_{L^2} + \| \rd_{A_0}\a_1 \|_{W^{-1,p}} \bigr)
\bigl( \| \a_2 \|_{L^2} + \| \rd_{A_0}\a_2 \|_{W^{-1,p}} \bigr) .
\end{align*}
This estimate continues to hold with a uniform constant if $A_0$ is replaced by an 
$L^p$-close connection $A$, since
\begin{equation}\label{A A0}
\|\rd_{A}\a-\rd_{A_0}\a\|_{W^{-1,p}}
\leq C \| [(A-A_0)\wedge\a] \|_{L^{\frac{2p}{p+2}}}
\leq C \| A-A_0  \|_{L^p} \| \a \|_{L^2}.
\end{equation}
Furthermore, the norm $\|\x\|_{W^{1,2}}$ is equivalent to
$ \|\x\|_{L^2} + \| \rd_{A_0} \x \|_{L^2}$ with a constant that depends continuously
on $A_0\in\cA^{0,p}(\Sigma)$ since 
$$
\|\rd_{A}\x-\rd_{A_0}\x\|_{L^2} \leq \| A-A_0  \|_{L^p} \| \x \|_{L^{\frac{2p}{p-2}}}
\leq  C \| A-A_0  \|_{L^p} \| \x \|_{W^{1,2}} .
$$
Finally, denote $q=\frac{2p}{p+2}$ then the dual exponents satisfy
$p'^{-1}-\frac 12 = q'^{-1}$, so we have the Sobolev embedding 
$W^{1,p'}\hookrightarrow L^{q'}$ and its dual $L^q\hookrightarrow W^{-1,p}$,
and thus with another uniform constant
$$
\|\rd_{A_0}\a_i \|_{W^{-1,p}}
\leq C \|\rd_{A_0}\a_i\|_{L^{\frac{2p}{p+2}}} .
$$
This proves (i), and we have moreover seen that
\begin{align*}
\| [ \a_1 \wedge \a_2 ] \|_{W^{-1,2}}\leq
C_3 \bigl( \| \a_1 \|_{L^2} + \| \rd_{A_0}\a_1 \|_{W^{-1,p}} \bigr) 
\bigl( \| \a_2 \|_{L^2} + \| \rd_{A_0}\a_2 \|_{W^{-1,p}} \bigr) .
\end{align*}
This estimate is not preserved under the gauge group, but it holds for every
smooth flat connection with a constant $C$ that depends on $A_0$.
If we consider any other connection $A\in\cA^{0,p}(\S)$, then the estimate continues
to hold with a new constant depending on $\|A-A_0\|_{L^p}$ by \eqref{A A0}.
\QED

\begin{cor} \label{cor div curl G} 
For every connection $A\in\cA^{0,p}(\S)$ there is a constant $C$
such that the following holds:
\begin{enumerate}
\item 
For any $\a\in\ker\rd_A\subset\Om^1(\S;\cg)$ and $\x\in\Om^0(\S,\cg)$
\begin{align*}
\| [ \a \wedge \rd_A\x ] \|_{W^{-1,2}}
&\leq
C \| \a \|_{L^2}  \| \rd_{A}\x \|_{L^2} .
\end{align*}
\item
For any $\x\in\Om^0(\S,\cg)$, $\a\in\Om^1(\S;\cg)$, 
and $u\in\cG(\S)$
$$
\left| \int_\S \tr\bigl( \x \cdot  \a \wedge \rd_{A} u \bigr) \right|
\;\leq\; 
C \| \x \|_{W^{1,2}} \bigl( \| \a \|_{L^2} + \| \rd_{A} \a \|_{W^{-1,p}} \bigr)
\| \rd_{A} u \|_{L^2} .
$$
\end{enumerate}
\end{cor}
\Pr
For (i) we assume without loss of generality that
$\x$ lies in some complement of $\ker\rd_A$, so we have
$\|\x\|_{L^q}\leq C_q \|\rd_A \x\|_{L^2}$ for any $1<q<\infty$ 
(but with a constant that might not depend continuously on $A$).
Now we apply Corollary~\ref{cor div curl}~(ii) with some $2<r<p$
to $\a\in\ker\rd_A\subset\Om^1(\S;\cg)$ and 
$\x\in\Om^0(\S,\cg)$ in the complement:
\begin{align*}
\| [ \a \wedge \rd_{A}\x ] \|_{W^{-1,2}} 
&\leq
C_1 \| \a \|_{L^2} 
\bigl( \| \rd_A\x \|_{L^2} + \| [ F_{A} , \x ] \|_{W^{-1,r}} \bigr) \\
&\leq
C_2 \| \a \|_{L^2} \| \rd_{A}\x \|_{L^2} .
\end{align*}
Here we estimated the curvature term as follows:
We used the Sobolev embedding 
$W^{1,r'}\hookrightarrow W^{1,p'}\hookrightarrow L^{\frac{2p'}{2-p'}}$
together with the H\"older inequality for
$L^p \cdot L^{\frac{2p'}{2-p'}}\hookrightarrow L^2$.
Moreover, we chose $q>1$ such that $\frac 1p+\frac 1q=\frac 1r$ and hence
\begin{align*}
\| [ F_{A} , \x ] \|_{W^{-1,r}}
&= \sup_{\|\psi\|_{W^{1,r'}}=1}  \biggl| \int_\S \la \rd_A\x \wedge \rd_A\psi \ra \biggr| \\
&\leq \sup_{\|\psi\|_{W^{1,r'}}=1} \Bigl( 
\| \rd_A\x \|_{L^2} \|[A,\psi]\|_{L^2} 
+ \| [A,\x] \|_{L^r } \|\rd\psi\|_{L^{r'}}  \Bigr) \\
&\leq 
C_3 \| \rd_A\x \|_{L^2} \|A\|_{L^p} 
+ \| A \|_{L^p } \| \x \|_{L^q} \\
&\leq C_A \| \rd_A\x \|_{L^2} .
\end{align*}
To prove (ii) we start with a smooth flat connection $A_0\in\cA_{\rm flat}(\S)$.
In that case we can use the Hodge decomposition $\a=\rd_{A_0}\z + \g$ with 
$\|\rd_{A_0}\z_i\|_{L^2}\leq \|\a_i\|_{L^2}$ and $\g\in\ker\rd_{A_0}^*$ such that
$$
\|\g\|_{L^p} \leq C_1 ( \|\a-\rd_{A_0}\z\|_{L^2} + \|\rd_{A_0}\g\|_{W^{-1,p}} )
\leq 2 C_1 ( \|\a\|_{L^2} + \|\rd_{A_0}\a\|_{W^{-1,p}} ) .
$$
Now use the H\"older inequality for $L^p\cdot L^2 \hookrightarrow L^{\frac{2p}{p+2}}$,
the Sobolev embedding $W^{1,2}(\S) \hookrightarrow L^{\frac{2p}{p-2}}(\S)$,
and Lemma~\ref{lem div curl} to estimate 
\begin{align*}
&\left| \int_\S \tr\bigl( \x \cdot  \a \wedge \rd_{A_0} u \bigr) \right| \\
& \leq \left| \int_\S \tr\bigl( \x \cdot  \rd_{A_0}\z \wedge \rd_{A_0} u \bigr) \right|
+ \left| \int_\S \tr\bigl( \x \cdot  \g \wedge \rd_{A_0} u \bigr) \right| \\
&\leq 
C_2 \|\x\|_{W^{1,2}} \| \rd_{A_0}\z\|_{L^2} \|\rd_{A_0} u \|_{L^2}
+ \|\x\|_{L^{\frac{2p}{p-2}}} \| \g \|_{L^p} \|\rd_{A_0} u\|_{L^2}  \\
&\leq 
C_3 \|\x\|_{W^{1,2}}  \bigl( \| \a \|_{L^2} + \| \rd_{A_0} \a \|_{W^{-1,p}} \bigr)
\| \rd_{A_0} u \|_{L^2} .
\end{align*}
Next, for general $A\in\cA^{0,p}(\S)$ we can assume without loss of generality that
$u$ lies in some complement of $\ker\rd_A$. Then we have
$\|u\|_{L^r}\leq C_4 \|\rd_A u\|_{L^2}$ for any fixed $1<r<\infty$.
Now 
we can estimate with $\frac 2r = 1 - \frac 12 + \frac 1p $
\begin{align*}
\left| \int_\S \tr\bigl( \x \cdot  \a \wedge \rd_{A} u \bigr) \right| 
&\leq 
\left| \int_\S \tr\bigl( \x \cdot  \a \wedge \rd_{A_0} u \bigr) \right|
+ \left| \int_\S \tr\bigl( \x \cdot  \a \wedge [(A-A_0), u] \bigr) \right| \\
&\leq 
C_3 \|\x\|_{W^{1,2}}  \bigl( \| \a \|_{L^2} + \| \rd_{A_0} \a \|_{W^{-1,p}} \bigr)
 \| \rd_{A_0} u \|_{L^2}  \\
&\quad
+ \|\x\|_{L^r} \|\a\|_{L^2} \|A-A_0\|_{L^p} \|u\|_{L^r} \\
&\leq 
C \|\x\|_{W^{1,2}}  \bigl( \| \a \|_{L^2} + \| \rd_{A} \a \|_{W^{-1,p}} \bigr)
\| \rd_{A} u \|_{L^2} .
\end{align*}

\eQED

%
%
%

The second crucial estimate for the local slice Theorem~\ref{thm local slice}
is the following.

\begin{lem} \label{lem est}
For every $A_0\in\cA^{0,p}(\S)$ there are constants $C$ and $\d>0$ so that the following
holds for all $A\in S_{A_0}$ in the local slice with $\|A-A_0\|_{L^2}\leq\d$ : 

If $\|u^*A-A_0\|_{L^2}\leq\d$ for some gauge transformation $u\in\cG^{1,p}(\S)$, then
$$
\| u^*A_0 - A_0 \|_{L^p} + \| A - A_0 \|_{L^p} \leq C \| u^*A - A_0 \|_{L^p} 
$$
and
$$
\| u^*A_0 - A_0 \|_{L^2} + \| A - A_0 \|_{L^2} \leq C \| u^*A - A_0 \|_{L^2} . 
$$
If moreover $u^*A\in S_{A_0}$ lies in the local slice
then automatically $u\in{\rm Stab}(A_0)$.
\end{lem}
\Pr
Let us write $A=A_0+a$ and $u^*A=A_0+b$, then we have
$$
\rd_{A_0}u = u(u^*A_0 - A_0)=
u ( u^* A  - A_0 - u^{-1}(A-A_0) u) = u b - a u
$$
with $\|a\|_{L^2}, \|b\|_{L^2}\leq\d$. 
Due to $\rd_{A_0}^*a=0$ we have
$$
\rd_{A_0}^*\rd_{A_0} u
= \rd_{A_0}^*(ub) - * ( * a\wedge\rd_{A_0}u ) .
$$
Now we can use the fact that $\rd_{A_0}^*:L^p\supset\im\rd_{A_0} \to W^{-1,p}$
is injective to estimate
\begin{align*}
\| \rd_{A_0} u \|_{L^p}
&\leq C_0 \| \rd_{A_0}^*\rd_{A_0} u \|_{W^{-1,p}} \\
&\leq C_0 \| \rd_{A_0}^*(ub) \|_{W^{-1,p}} + C_0 \| * a\wedge\rd_{A_0}u  \|_{W^{-1,p}} \\
&\leq C_1 \| u b \|_{L^p} + C_1 \| a \|_{L^2} \| \rd_{A_0}u  \|_{L^p} .
\end{align*}
Here we used the Sobolev embedding 
$W^{1,p'}\hookrightarrow L^{q'}$ with $\frac 1{q'}=\frac 1{p'}-\frac 12$
and its dual $L^q\hookrightarrow W^{-1,p}$ with $\frac 1q=\frac 1p + \frac 12$.
If we pick $\d<C_1$ then this proves the first part of the first inequality,
$$
(1-C_1\d) \| u^*A_0 - A_0 \|_{L^p} \leq C_1 \| b \|_{L^p} = C_1 \| u^*A - A_0 \|_{L^p},
$$
For the second part just use $au = ub - \rd_{A_0}u$ to see that
$$
\|A-A_0\|_{L^p} = \|a u\|_{L^p} \leq \| b\|_{L^p} + \|\rd_{A_0} u\|_{L^p} .
$$
The second inequality, with $p$ replaced by $2$, is proven in the same way.
The crucial estimate now uses Corollary~\ref{cor div curl G}~(ii)
with $\rd_{A_0}^*a=0$,
\begin{align*}
\| \rd_{A_0}^*\rd_{A_0} u \|_{W^{-1,2}}
&= \sup_{\x\neq 0} \; \|\x\|_{W^{1,2}}^{-1} 
\biggl| \int_\S \la \x ,
\bigl( \rd_{A_0}*(ub) +  * a\wedge\rd_{A_0}u  \bigr) \ra \biggl| \\
&\leq C_3 \| u b \|_{L^2} + C_3 \| a \|_{L^2} \| \rd_{A_0}u  \|_{L^2} .
\end{align*}
If we moreover assume $\rd_{A_0}^*b=0$ then
$$
\rd_{A_0}^*\rd_{A_0} u
= - * ( \rd_{A_0}u \wedge * b) - * ( * a\wedge\rd_{A_0}u ) 
$$
and we can estimate as before
\begin{align*}
\| \rd_{A_0} u \|_{L^p}
&\leq C_0 \bigl(  \| \rd_{A_0}u \wedge *b \|_{W^{-1,p}}
+ \| * a\wedge\rd_{A_0}u  \|_{W^{-1,p}} \bigr) \\
&\leq C_1 \bigl(\| a \|_{L^2} + \| b \|_{L^2} \bigr) \| \rd_{A_0}u  \|_{L^p}
\;\leq\; 2 C_1 \d \| \rd_{A_0}u  \|_{L^p} .
\end{align*}
For $\d<\half C_1^{-1}$ this implies $\rd_{A_0}u=0$, i.e.\ $u^*A_0=A_0$.
\QED

\section{$L^2$-topology on the space of connections}
\label{proofs}

As before, let $\S$ be a closed Riemannian surface, let $\rG$ be a compact Lie group with Lie algebra $\cg$, and consider a trivialized $\rG$-bundle over $\S$.
The aim of this section is to provide some understanding of the $L^2$-topology
on the space of connections $\cA^{0,p}(\S)=L^p(\S;T^*\S\otimes\cg)$ and its quotient by the
gauge group $\cG^{1,p}(\S)=W^{1,p}(\S,\rG)$ for fixed $p>2$. 
Firstly, we will prove the local slice Theorem~\ref{thm local slice} with the following refinement.

\begin{rmk} \label{rmk slice}
The differential of the map $\cm$ in Theorem~\ref{thm local slice} is bounded
in the $L^2$-topology in the following sense:
For any smooth path $A:(-1,1)\to B_\d(A_0)$ let 
$(A_0+a,u):(-1,1)\to S_{A_0}(\ep)\times\cG^{1,p}(\S)$ be a smooth representative
of $\cm^{-1}\circ A$. Then
$$
\|\rd_{A_0}( u^{-1} \pd_t u ) \|_{L^2} 
+ \|\pd_t a - [a, u^{-1} \pd_t u ] \|_{L^2} 
\leq C \|\pd_t A\|_{L^2} .
$$
\end{rmk}

Based on this local slice theorem, we will prove the following quantitative local pathwise connectedness and local quasiconvexity of the gauge orbits in the $L^2$-topology.

\begin{thm} \label{thm short}
For every connection $A_0\in\cA^{0,p}(\S)$ there exist constants $C$ and $\d>0$ 
so that the following holds:
For any $A\in S_{A_0}$ and $u\in\cG^{1,p}(\S)$ 
with $\|A-A_0\|_{L^2}<\d$ and $\|u^*A-A_0\|_{L^2}<\d$ there exists a smooth path
$v:[0,1]\to\cG^{1,p}(\S)$ with $v(0)={\emph\one}$ such that
$v(1)^{-1\;*}u^*A \in S_{A_0}$ and 
\begin{equation} \label{short}
 \bigl\| \pd_t \bigl( v(t)^*A \bigr) \bigr\|_{L^2}
\;\leq\; C \bigl( 1+ \|A-A_0\|_{L^p} \bigr) \|u^*A-A\|_{L^2}
\qquad \forall t\in[0,1] .
\end{equation}
Here $u \cdot v(1)^{-1}\in {\rm Stab}(A_0)$, so if $A=A_0$ or if $A_0$ is irreducible, then moreover $v(1)^*A=u^*A$, and hence
$t\mapsto v(t)^*A$ connects $A$ to $u^*A$ by a smooth path in the gauge orbit of 
$L^2$-length bounded by $C(1+\|A-A_0\|_{L^p})\|u^*A-A\|_{L^2}$.
\end{thm}


\begin{rmk}
We allow for $A\neq A_0$ in Theorem~\ref{thm short} in order to obtain local pathwise connectedness statements with uniform constants on $L^p$-balls.
Near a reducible connection $A_0$ however, our method only provides a path in the gauge orbit from $A$ to $v(1)^*A$ that might differ from $u^*A$ by an element of the stabilizer
${\rm Stab}(A_0)$.
So the question of local connectedness with uniform constants
reduces to the action of ${\rm Stab}(A_0)$ on the local slice:

We have $w:=u \cdot v(1)^{-1}\in {\rm Stab}(A_0)$ such that $A, w^*A\in S_{A_0}$
are $L^2$-close.
Can these be connected within the gauge orbit in $S_{A_0}$ 
with a constant in (\ref{short}) that only depends on $A_0$?
The difficulty in this question is to obtain uniform bounds for the action of 
the compact group ${\rm Stab}(A_0)$ on the noncompact local slice $S_{A_0}$.
In our application to gauge invariant Lagrangian submanifolds $\cL\subset\cA(\S)$ 
we can bypass this difficulty by establishing that $\cL\cap S_{A_0}$ is a 
finite dimensional manifold near $A_0$. 
\end{rmk}

Before going into the more technical proofs, let us explain how the local connectivity and quasiconvexity in Theorem~\ref{thm im Kleinen} follows from the special case of $A_0=A$ of
Theorem~\ref{thm short}.

\medskip
\noindent{\bf Proof of Theorem~\ref{thm im Kleinen} :}
A simple corollary of  Theorem~\ref{thm short} is the following: 
Given $A_0\in\cA^{0,p}(\S)$ there exist constants $C$ and $\delta>0$ such that for any $u^*A_0 \in\cG^{1,p}(\S)^*A_0$  with $\|u^*A_0 - A_0\|_{L^2}\leq\delta$ there exists a smooth path $[0,1]\to \cG^{1,p}(\S)^*A_0\subset\cA^{0,p}(\S)$, $t\mapsto A_t=v(t)^*A_0$ 
connecting $A_0$ to $A_1=u^*A_0$ 
with derivative $\|\partial_t A_t\|_{L^2} \leq C \|u^*A_0-A_0\|_{L^2}$. 

This proves Theorem~\ref{thm im Kleinen}~(i) and (ii) in the case $A_0=B$; the general case follows from gauge invariance as follows:
For an arbitrary base point $B\in\cA^{0,p}(\S)$, let $C$ and $\delta$ be the above constants for $A_0=B$. If we now consider any $A_0=u_0^*B$ and $A_1=u_1^*B$ with $\|u_1^*B - u_0^*B\|_{L^2}\leq\delta$, then
$\|(u_1u_0^{-1})^*B -B\|_{L^2}=\|u_1^*B -u_0^*B\|_{L^2}\leq\delta$ and we obtain a smooth path $B_t$ connecting $B_0=B$ to $B_1=(u_1u_0^{-1})^*B$. We can transform it back to a path $A_t=u_0^*B_t$
connecting $A_0$ and $A_1$ with derivative 
$$
\|\pd_t (u_0^* B_t) \|_{L^2}
= \|u_0^{-1} ( \pd_t B_t ) u_0 \|_{L^2}
= \| \pd_t B_t  \|_{L^2}
\leq C \|(u_1u_0^{-1})^*B -B\|_{L^2} 
= C \|A_1-A_0\|_{L^2} .
$$
Finally, to prove (iii) suppose that $B_0\in\cA^{0,p}(\S)$ is irreducible, then for any $L^2$-close base point $B$ with $\|B-B_0\|_{L^2}\leq\delta$ we find paths with derivative control
$\|\partial_t A_t\|_{L^2} \leq C (1+ \|B-B_0\|_{L^p}) \|u^*A-A\|_{L^2}$ for any $A\in\cG^{1,p}(\S)^*B$. Clearly, the constant 
$C (1+ \|B-B_0\|_{L^p})$ can be replaced by a uniform constant 
for all $B$ in an $L^p$-bounded subset of the given $L^2$-neighbourhood of $B_0$.
\QED

In the remainder of this section, we give proofs of the local slice theorem and quantitative local quasiconvexity,
using estimates from Section \ref{est} which are based on a version of the div-curl lemma.

\medskip
\noindent{\bf Proof of Theorem~\ref{thm local slice} : }
The map $\cm$ in \eqref{map} is well defined. In particular, we have 
$\cm[(g^*(A_0+a), u g)]=\cm[(A_0+a, u)]$ for all $g\in{\rm Stab}(A_0)$.
It is equivariant in the sense that
$\cm[(A_0+a, v^{-1} u)]=v^*\,\cm[(A_0+a, u)]$ for all $v\in\cG^{1,p}(\S)$.

In a first step we will show that the differential $D\cm_{(A_0+a,u)}$ of $\cm$
is injective for all $(a,u)\in \cS_{A_0}(\ep)\times\cG^{1,p}(\Sigma)$ with $\ep>0$ sufficiently small.
By equivariance,
$D\cm_{(A_0+a,u)}(\alpha,\xi) = u \bigl(D\cm_{(A_0+a,\smone)}(\alpha,\xi)\bigr) u^{-1}$
it suffices to prove injectivity for $u=\one$, where
$$
D\cm_{(A_0+a,\smone)} : (\a,\x) \mapsto \a - \rd_{A_0}\x - [a,\x] 
$$ 
acts on $(\a,\x)\in \ker(\rd_{A_0}^*)\times W^{1,p}(\Sigma,\cg)$
in the $L^2$-orthogonal complement of ${\rm Stab}(A_0)^*(A_0+a,\one)$, i.e.\ 
$(\a,\x)\perp\{ ([a,\psi],\psi ) \st \psi\in\ker\rd_{A_0} \}$.
So let $(\a,\x)\in\ker D\cm_{(A_0+a,\smone)}$, then $\rd_{A_0}\x = \a - [a,\x]$ and hence
\begin{align*}
\rd_{A_0}^* ( \rd_{A_0}\x)  
=   \rd_{A_0}^* \a - \bigl[ \rd_{A_0}^*a , \x \bigr] - * \bigl[ *a \wedge \rd_{A_0}\x \bigr]
=  - *\bigl[*a\wedge\rd_{A_0}\x \bigr] .
\end{align*}
Here we can use the H\"{o}lder inequality for $\frac 1p + \frac 1r = \frac 12$
and the Sobolev embedding $W^{1,p'}(\S)\hookrightarrow L^r(\S)$ due to 
$\frac 1{p'} - \frac 12 = \frac 1r$ to estimate with the Sobolev constant $C_1$
\begin{align*}
\| *[*a\wedge\rd_{A_0}\x] \|_{W^{-1,p}}
&=\sup_{\psi\neq 0} \bigl| \tint_\S \la *[*a\wedge\rd_{A_0}\x] \,,\, \psi \ra \bigr|
\|\psi\|_{W^{1,p'}}^{-1} \\
&\leq \sup_{\psi\neq 0} \| a \|_{L^2} \|\rd_{A_0}\x\|_{L^p} \| \psi \|_{L^r} \|\psi\|_{W^{1,p'}}^{-1} 
\leq C_1\|a\|_{L^2} \|\rd_{A_0}\x\|_{L^p} .
\end{align*}
With Lemma~\ref{lem:c0} this implies
$$
\|\rd_{A_0}\x\|_{L^p}
\leq C_0 \|\rd_{A_0}^* ( \rd_{A_0}\x) \|_{W^{-1,p}}
\leq C_0 C_1\|a\|_{L^2} \|\rd_{A_0}\x\|_{L^p} .
$$
If we assume $\|a\|_{L^2} < (C_0 C_1)^{-1}$ then we can conclude $\rd_{A_0}\x=0$, and hence $\a=[a,\xi]$.
Since $(\alpha,\xi)\perp([a,\xi],\xi)$ this implies $(\a,\xi)=0$.
This proves the injectivity of $D\cm_{(A_0+a,u)}$ 
for all $(a,u)\in \cS_{A_0}(\ep)\times\cG^{1,p}(\Sigma)$ with $0<\ep<(C_0 C_1)^{-1}$.

Secondly, the differentials $D\cm_{(A_0+a,u)}$ are Fredholm maps that 
vary smoothly with $(A_0+a,u)\in S_{A_0}\times\cG^{1,p}(\S)$.
The differential at $(A_0,\one)$, given by $D\cm_{(A_0,\smone)}(\a,\xi)=\a-\rd_{A_0}\xi$, is surjective
by the Hodge decomposition $L^p(\S,\rT^*\S\otimes\cg)=\ker\rd_{A_0}^* \oplus \rd_{A_0} W^{1,p}(\S,\cg)$.
So, by the stability of the Fredholm index, all differentials for $(a,u)\in \cS_{A_0}(\ep)\times\cG^{1,p}(\Sigma)$
are bijections, and hence $\cm$ is a local diffeomorphism onto its image.
Indeed, $\cm$ is also injective if we choose $0<\ep<\delta$ with the $\delta>0$ from Lemma~\ref{lem est}:
If $u^{-1\;*}(A_0+a)=\tilde{u}^{-1\;*}(A_0+\tilde{a})$ then 
$(u^{-1}\tilde u)^*(A_0+a)=A_0+\tilde{a}\in S_{A_0}$, and hence the Lemma implies 
$g:=u^{-1}\tilde u\in{\rm Stab(A_0)}$. So we have $(A_0+\tilde{a},\tilde{u}) = g (A_0 + a , u )$, i.e.\ 
the two pairs are equivalent in $\bigl(\cS_{A_0}(\ep)\times\cG^{1,p}(\Sigma)\bigr)/{\rm Stab}(A_0)$.
This proves that $\cm$ is a diffeomorphism onto its image.

Finally, our aim is to prove that $B_\d:=\{ A \in\cA^{0,p}(\S) \st \|A-A_0\|_{L^2}<\d \}$
is a subset of $\im\cm$ for appropriate choices of $\ep$ and $\d$.
This will follow from a connectedness argument: 
Note that $B_\d$ is connected and $B_\d\cap\im\cm$ is nonempty since it contains 
$A_0=\cm(A_0,\one)$. 
So if $B_\d\cap\im\cm$ is both open and closed with respect to the $L^p$-topology 
on $B_\d$, then it has to be the whole space $B_\d$, as claimed.

The intersection is open since both 
$\im\cm$ and $B_\d$ are open subsets of $\cA^{0,p}(\S)$.
To see that $B_\d\cap\im\cm$ is closed in $B_\d$ consider an $L^p$-convergent sequence
$B_\d\cap\im\cm\ni A_i=u_i^*(\hat A_i)\rightarrow A_\infty\in B_\d$ 
with $\hat A_i\in \cS_{A_0}(\ep)$.
If we choose $\ep>0$ sufficiently small then Lemma~\ref{lem est} applies to give
$$
\| u_i^*A_0 - A_0 \|_{L^p} + \| \hat A_i - A_0 \|_{L^p} 
\leq C \| A_i - A_0 \|_{L^p} .
$$
Since the right hand side is bounded we find weakly convergent subsequences
$\hat A_i \rightharpoonup \hat A_\infty\in\cA^{0,p}(\S)$ and 
$u_i \rightharpoonup u_\infty\in\cG^{1,p}(\S)$ (with strong $\cC^0$-convergence). 
The weak convergence preserves the local slice condition 
$\rd_{A_0}^*(\hat A_\infty - A_0) = 0$ and the 
unique limit $u_\infty^*\hat A_\infty = A_\infty$.
Moreover, again from Lemma~\ref{lem est},
$$
\|\hat A_\infty - A_0 \|_{L^2} \leq
 \liminf_{i\to\infty}  \| \hat A_i - A_0 \|_{L^2} \leq
 \lim_{i\to\infty} C \| A_i - A_0 \|_{L^2} \leq C\d .
$$
So if we choose $\d<C^{-1}\ep$, then $\hat A_\infty$ automatically lies in $\cS_{A_0}(\ep)$
and hence we have $A_\infty=(u_\infty)^{-1\;*}\hat A_\infty \in\im\cm$.
This proves the closedness of  $B_\d\cap\im\cm\subset B_\d$ and thus finishes the proof.
\QED

\noindent{\bf Proof of Remark \ref{rmk slice} : }
We are considering paths $u:(-1,1)\to\cG^{1,p}(\S)$, $A_0+a:(-1,1)\to S_{A_0}(\ep)$,
and $A:(-1,1)\to\cA^{0,p}(\S)$ such that and  $u^*A=A_0+a$. Differentiating this, we obtain
$\rd_{u^*A}(u^{-1}\pd_t u) + u^{-1}\pd_t A \, u = \pd_t a $
and hence 
\begin{align*}
\rd_{A_0}(u^{-1}\pd_t u) = \pd_t a - u^{-1}\pd_t A \, u - [ a , u^{-1}\pd_t u ] .
\end{align*}
From this we calculate, using the Coulomb gauge $\rd_{A_0}^*a = - *\rd_{A_0} *a = 0$,
\begin{equation*}
\rd_{A_0}^*\rd_{A_0} (u^{-1} \pd_t u ) = 
- \rd_{A_0}^*( u^{-1} \pd_t A \, u) 
- * [ *a \wedge \rd_{A_0}( u^{-1}\pd_t u ) ] .
\end{equation*}
Now we use Lemma~\ref{lem:c0} and Corollary~\ref{cor div curl G}~(i) and 
to obtain (with another constant $C_2$)
\begin{align*}
& \| \rd_{A_0}( u^{-1} \pd_t u ) \|_{L^2} \\
&\leq C_0 \| \rd_{A_0}^*\rd_{A_0}( u^{-1} \pd_t u ) \|_{W^{-1,2}} \\
&\leq C_0 \bigl( \| \rd_{A_0}^*( u^{-1} \pd_t A u) \|_{W^{-1,2}}
+ \bigl\| \bigl[ * a \wedge \rd_{A_0} (u^{-1}\pd_t u) \bigr]\bigr\|_{W^{-1,2}} \bigr)\\
&\leq C_2 \| u \pd_t A  u^{-1} \|_{L^2}
+ C_2 \| * a \|_{L^2}  \| \rd_{A_0} (u^{-1}\pd_t u) \|_{L^2} \\
&\leq C_2 \| \pd_t A \|_{L^2}
+ C_2 \ep  \| \rd_{A_0} (u^{-1}\pd_t u) \|_{L^2} .
\end{align*}
Here we have $a = u^*A-A_0\in S_{A_0}(\ep)$ and we can choose this $L^2$-ball in the 
local slice sufficiently small, $\ep\leq \half C_2^{-1}$, to obtain
$\| \rd_{A_0}( u^{-1} \pd_t u ) \|_{L^2}
\;\leq\; 2 C_2 \| \pd_t A \|_{L^2}$.

Next, recall that we have
$\pd_t a - [a,u^{-1}\pd_t u] = u^{-1} \pd_t A \, u + \rd_{A_0} (u^{-1}\pd_t u) $.
Taking the $L^2$-norm on both sides and using the previous estimate now gives
$$
\|\rd_{A_0}( u^{-1} \pd_t u ) \|_{L^2} 
+ \|\pd_t a - [a, u^{-1} \pd_t u ] \|_{L^2} 
\leq (1+4C_2) \|\pd_t A\|_{L^2} .
$$

\eQED

\noindent{\bf Proof of Theorem \ref{thm short} : }
The smooth path $B:[0,1]\to\cA^{0,p}(\S)$ given by
$B(t)=A + t(u^*A - A)$ obviously connects $B(0)=A$ to $B(1)=u^*A$
within $B_\delta(A_0)$ and has $L^2$-speed $\|u^*A - A\|_{L^2}$.
The idea is to 
construct the path $v:[0,1]\to\cG^{1,p}(\S)$
by finding gauge transformations $v(t)^{-1}$ that take $B(t)$
into the local slice at $A_0$. Then we will be able to use Remark~\ref{rmk slice} to control the length of the path $v^*A$.

More precisely, we pick $\delta>0$ as in Theorem~\ref{thm local slice} to obtain 
a smooth path $\g:[0,1]\to(\cS_{A_0}(\ep)\times\cG^{1,p}(\S))/{\rm Stab}(A_0)$ with
$\cm(\g(t))=B(t)$ and $\g(0)=[(0,\one)]$.
We can project it to $\cG^{1,p}(\S)/{\rm Stab}(A_0)$ and
then lift it to a smooth path $w:[0,1]\to\cG^{1,p}(\S)$
starting at $w(0)=\one$ and solving $\rd_{A_0}^*(w^*B - A_0)=0$.
Now $v(t):= w(t)^{-1}$ also defines a smooth path $v:[0,1]\to\cG^{1,p}(\S)$
starting at $v(0)=\one$, which for $t=1$ clearly satisfies
$v(1)^{-1\;*} B(1) =  v(1)^{-1\;*} u^*A \in S_{A_0}$.
For reducible $A_0$ we can moreover assume that
$w^{-1}\pd_t w$ is $L^2$-orthogonal to $\ker\rd_{A_0}$ 
after the following modification of $w$:

Let $\pi:W^{1,p}(\S,\cg)\to W^{1,p}(\S,\cg)$ be the $L^2$-orthogonal projection
to $\ker\rd_{A_0}$. Then we solve $\pd_t g \cdot g^{-1}=-\pi(w^{-1}\pd_t w)$
by a smooth path $g:[0,1]\to{\rm Stab}(A_0)\subset \cG^{1,p}(\S)$ with $g(0)=\one$.
(Note that $g$ automatically takes values in the stabilizer since
$\pd_t(g^*A_0)=g^{-1}\rd_{A_0}(\pd_t g\cdot g^{-1}) g = 0$.)
Now $\tilde w:=wg:[0,1]\to\cG^{1,p}(\S)$ satisfies
$\tilde w^{-1}\pd_t \tilde w = g^{-1}\bigl( w^{-1}\pd_t w  - \pi(w^{-1}\pd_t w) \bigr) g$,
and this is orthogonal to $\ker\rd_{A_0}$ since the latter is invariant under conjugation with $g\in{\rm Stab}(A_0)$.
In addition, $\tilde w^{-1\;*}A_0 = w^{-1\;*}g^{-1\;*}A_0 = w^{-1\;*}A_0$ still satisfies the Coulomb gauge condition.

Now in order to control the length of the path $t\mapsto v(t)^*A$, first notice that
$$
 \pd_t ( v^*A ) = \rd_{v^*A}( v^{-1} \pd_t v ) 
 = v^{-1} \bigl( \rd_{A}( \pd_t v \, v^{-1} ) \bigr) v
 = - v^{-1} \bigl( \rd_{A}( w^{-1} \pd_t w ) \bigr) v .
$$
From Remark~\ref{rmk slice} we have 
$\| \rd_{A_0}( w^{-1} \pd_t w ) \|_{L^2} \leq C  \| \pd_t B \|_{L^2} $, and hence
\begin{align*}
\| \pd_t ( v^*A ) \|_{L^2} 
&\leq \| \rd_{A_0}( w^{-1} \pd_t w ) \|_{L^2} 
+ \| A-A_0 \|_{L^p} \| w^{-1} \pd_t w \|_{L^{\frac{2p}{p-2}}} \\
&\leq \bigl( 1+ C_1\| A-A_0 \|_{L^p} \bigr)
\| \rd_{A_0}( w^{-1} \pd_t w ) \|_{L^2} \\
&\leq C \bigl( 1+ \| A-A_0 \|_{L^p} \bigr)  \| \pd_t B \|_{L^2} 
\;=\; C \bigl( 1+ \| A-A_0 \|_{L^p} \bigr) \|u^*A - A\|_{L^2} .
\end{align*}
Here we also used 
$\| \x \|_{L^{\frac{2p}{p-2}}}\leq C_1 \| \rd_{A_0} \x \|_{L^2}$
for $\x=w^{-1}\pd_t w\in(\ker\rd_{A_0})^\perp$.

Finally, we have $(u\cdot v(1)^{-1})^*A = w(1)^*B(1) \in S_{A_0}(\ep)$ by construction.
Since $A$ also lies in the local slice and $L^2$-close to $A_0$, 
the local uniqueness in Lemma~\ref{lem est} for sufficiently small $\ep,\delta>0$ 
implies that $u \cdot v(1)^{-1} \in{\rm Stab}(A_0)$.
So if $A_0$ is irreducible or if $A=A_0$ then
$(u\cdot v(1)^{-1})^*A = A$ and hence $u^*A = v(1)^*A$.
\QED

\section{Gauge invariant Lagrangians in the space of connections}
\label{sec Lag}

In this section we consider gauge invariant Lagrangian submanifolds in
the space of connections $\cA(\S)$ over a closed Riemann surface.
These are discussed in detail in \cite[Section~4]{W Cauchy}, where
we established their basic structure in the $L^p$-topology for $p>2$.
The aim of this section is to provide some control of their $L^2$-geometry
despite the fact that it is unclear whether their $L^2$-closure is even a topological manifold.
We will work with the following general class of Lagrangians into which for example all handle body Lagrangians $\cL_H$ fall.

%
%

\begin{dfn}\label{dfn L}
We call $\cL\subset\cA^{0,p}(\Sigma)$ a {\bf gauge invariant Lagrangian submanifold}
if it is a Banach submanifold,
invariant under the action of $\cG^{1,p}(\Sigma)$, and Lagrangian 
in the following sense:
For every $A\in\cL$ the tangent space 
$\rT_A\cL\subset L^p(\Sigma,\rT^*\Sigma\otimes\cg)$
is Lagrangian, i.e.\ for every $\alpha\in L^p(\Sigma,\rT^*\Sigma\otimes\cg)$
\begin{equation}\label{lagr}
\o(\alpha,\beta):=\int_\Sigma \la \alpha \wedge \beta \ra = 0\quad
\forall \beta\in \rT_A\cL\qquad
\iff\qquad \alpha\in \rT_A\cL .
\end{equation}
We moreover assume\footnote{
This follows directly from the other assumptions if $G$ has discrete center, i.e.\ $[\cg,\cg]=\cg$; see \cite[Section~4]{W Cauchy}.
}
that $\cL$ lies in the subset of weakly flat connections, 
$\cL\subset\cA^{0,p}_{\mathrm{flat}}(\Sigma)$,
and that the quotient space $\cL/\cG^{1,p}_z(\S)$ by the based gauge group is compact. 
\end{dfn}

We know from \cite[Section~4]{W Cauchy} that any such Lagrangian $\cL$  
is a totally real submanifold with respect to the Hodge $*$ 
operator for any metric on $\Sigma$, i.e.\ for every $A\in\cL$
$$
L^p(\Sigma,\rT^*\Sigma\otimes\cg)=\rT_A\cL\oplus *\rT_A\cL .
$$ 
This should be compared with the Hodge decomposition for any 
$A\in\cA^{0,p}_{\rm flat}(\S)$,
$$
L^p(\Sigma,\rT^*\Sigma\otimes\cg)
= \im\rd_{A} \oplus h^1_{A} \oplus *\im\rd_{A} 
$$
with $h^1_A=\ker\rd_A\cap\ker\rd_A^*$.
The quotient $L:=\cL/\cG(\Sigma)$ has singularities in general,
but $\cL$ has the structure of a principal bundle
\begin{equation}\label{bundle}
\cG^{1,p}_z(\Sigma)\hookrightarrow \cL \to \cL/\cG^{1,p}_z(\Sigma)
\end{equation}
over a smooth quotient $\cL/\cG^{1,p}_z(\Sigma)$.
Here we fix a base point set $z\subset\Sigma$ consisting of exactly 
one point in each connected component of $\Sigma$, 
then the fibre is the based gauge group
$\cG^{1,p}_z(\Sigma)=\{u\in\cG^{1,p}(\Sigma)\,|\, u(z)\equiv\one\}$.

Every class in $\cL/\cG^{1,p}_z(\Sigma)$ has a smooth representative.
So the bundle structure shows that the $W^{k,q}$-closure or restriction
of $\cL$ is again a smooth Banach submanifold of $\cA^{k,q}(\S)$ as long
as $(k+1)q>2$, so $\cG^{k,q}_z(\S)$ is well defined.
It is however unclear whether the $L^2$-closure of $\cL$ will necessarily
be a smooth Hilbert submanifold of $\cA^{0,2}(\S)$.
The lack of $L^2$-charts for $\cL$ will be compensated by the following
proposition and the local slice Theorem~\ref{thm local slice}.
Recall the definition of the $L^2$-balls in the local slice at 
$A_0\in\cA^{0,p}(\S)$,
$$
S_{A_0}(\ep):= 
\bigl\{ A \in\cA^{0,p}(\S) \st \rd_{A_0}^*(A-A_0) = 0 , \|A-A_0\|_{L^2}<\ep \bigr\} .
$$

\begin{prp} \label{prp lag slice}
For any $A_0\in\cL$ there is an $\ep>0$ such that the intersection
$L_{A_0}:=\cL\cap S_{A_0}(\ep)$ of the Lagrangian with the $L^2$-ball 
in the local slice
is a submanifold of dimension $\dim L_{A_0}=\half\dim h^1_{A_0}$.
\end{prp}
\Pr
Since $\cL\subset\cA^{0,p}_{\rm flat}(\S)$ we have the Hodge decomposition
\cite[Lemma~4.1]{W Cauchy} in the $L^p$-topology,
$$
\Om^1(\S,\cg) = \ker\rd_{A_0}^* \oplus \im\rd_{A_0} .
$$
This shows that $\rT_{A_0} S_{A_0}=\ker\rd_{A_0}^*$ is transverse to
$\rT_{A_0}\cL\supset\im\rd_{A_0}$.
We claim that the transversality of $S_{A_0}$ and $\cL$ persists in an
$L^2$-neighbourhood of $A_0$.
Then by the implicit function theorem the intersection 
$\cL\cap S_{A_0}(\ep)$ is a submanifold of both $\cL$ and $S_{A_0}$.
Note that $\rT_A S_{A_0}=\ker\rd_{A_0}^*$ for all $A\in S_{A_0}$ and
$\rT_{A}\cL\supset\im\rd_{A}$ for all $A\in\cL$.
So it suffices to prove that for all $A\in\cL$ with $\rd_{A_0}^*(A-A_0)=0$ 
and $\|A-A_0\|_{L^2}\leq\ep$ sufficiently small we have
$$
\Om^1(\S,\cg) = \ker\rd_{A_0}^* + \im\rd_{A} .
$$
To prove this we need to consider any $\a\in\Om^1(\S,\cg)$ and 
find $\x\in W^{1,p}(\S,\cg)$ such that $\a-\rd_A\x\in\ker\rd_{A_0}^*$.
This is achieved by solving $\rd_{A_0}^*\rd_A\x=\rd_{A_0}^*\a$, 
so we only need to check the surjectivity of
$$
\rd_{A_0}^*\rd_A : W^{1,p}(\S,\cg)\supset(\ker\rd_{A_0})^\perp \to 
\im\rd_{A_0}^*\subset W^{-1,p}(\S,\cg) .
$$
For $A=A_0$ this is a Fredholm operator of index $0$.
For any other $A\in\cA^{0,p}(\S)$ it is a compact perturbation 
(and thus also Fredholm of index $0$) since
$$
\|\rd_{A_0}^*\rd_A\x - \rd_{A_0}^*\rd_{A_0}\x \|_{W^{-1,p}}
\leq C \| (A-A_0) \|_{L^p} \| \x \|_{L^\infty} 
$$
and $W^{1,p}(\S)\hookrightarrow L^\infty(\S)$ is compact.
So instead of the surjectivity we can check the injectivity:
Let $\x\in(\ker\rd_{A_0})^\perp\subset  W^{1,p}(\S,\cg)$
with $\rd_{A_0}^*\rd_A\x=0$, then we use Corollary~\ref{cor div curl}~(i)
with the local slice condition $\rd_{A_0}^*(A-A_0)=0$ and the weak
flatness $\rd_{A_0}\rd_{A_0}\xi = 0$ of $A_0$ to estimate
\begin{align*}
0 &= \int_\S \la \rd_{A_0}\x \wedge * \rd_A\x \ra \\
&= \| \rd_{A_0}\x \|_{L^2}^2 
+ \int_\S \la \rd_{A_0}\x \wedge * [A-A_0,\x] \ra \\
&= \| \rd_{A_0}\x \|_{L^2}^2 
+ \int_\S \la \x , [*(A-A_0)\wedge\rd_{A_0}\x] \ra \\
&\geq \| \rd_{A_0}\x \|_{L^2}^2 
- C_1 \bigl(\|\x\|_{L^2} + \|\rd_{A_0}\x\|_{L^2} \bigr)
\|A-A_0\|_{L^2} \|\rd_{A_0}\x\|_{L^2} \\
&\geq \bigl( 1 - 2C_1 C_2 \|A-A_0\|_{L^2} \bigr) \|\rd_{A_0}\x\|_{L^2} .
\end{align*}
Here we used the estimate $\|\x\|_{L^2}\leq C_2\|\rd_{A_0}\x\|_{L^2}$
for $\x\in(\ker\rd_{A_0})^\perp$ and some $C_2\geq1$.
This calculation shows the injectivity of $\rd_{A_0}^*\rd_A$ 
and thus the claimed transversality for $\|A-A_0\|_{L^2}<\ep:=(2C_1C_2)^{-1}$.
\QED

The next two results replace the $L^2(\S)$ to $L^3(H)$ extension properties \cite[Lemma~1.6]{W bubb} of handle body Lagrangians $\cL_H\subset\cA(\S)$.
First, we have the following weak form of a uniform curvature bound, restating Lemma~\ref{lem lin ext i} from the introduction.

\begin{lem} \label{lem lin ext}
There is a constant $C_{\rT\cL}$
such that any smooth path $A:(-s_0,s_0)\to\cL$ 
satisfies 
\begin{align*}
\int_\S \la \pd_s A(0) \wedge \pd_s\pd_s A(0) \ra 
&\leq C_{\rT\cL} \bigl\| \pd_s A(0) \bigr\|_{L^2(\Si)}^3  .
\end{align*}
\end{lem}
\Pr
This estimate is preserved under constant gauge transformations in
$\cG^{1,p}(\S)$. So by the compactness of $\cL/\cG^{1,p}_z(\S)$ it
suffices to establish the estimate for paths $A:(-s_0,s_0)\to\cL$ that pass
through the local slice, $A(0)\in S_{A_0}$, 
for some fixed smooth $A_0\in\cL$.
We can moreover assume that the entire path $A$ 
lies in an $L^p$-neighbourhood of $A_0$.
Then in Theorem~\ref{thm local slice} we can replace
$\cG^{1,p}(\S)/{\rm Stab}(A_0)$ by the image of the exponential map on
a $W^{1,p}$-ball in the $L^2$-complement of $\ker\rd_{A_0}$, that is
$D^{1,p}(\d)\cap (\ker\rd_{A_0})^\perp \subset W^{1,p}(\S,\cg)$
with
$D^{1,p}(\d):=\{\xi\in W^{1,p}(\S,\cg) \st \|\xi\|_{W^{1,p}}<\d\}$.
For sufficiently small $\d>0$ the map
\[
\cm:\; \begin{aligned}
\cS_{A_0}(\ep) \times \bigl( D^{1,p}(\d)\cap(\ker\rd_{A_0})^\perp\bigr) 
\; & \to & \cA^{0,p}(\S) \\
(A_0+a,\x) \;& \mapsto & \exp(\x)^*(A_0 + a)
\end{aligned}
\]
is a diffeomorphism onto its image, which contains an $L^p$-neighbourhood
of $A_0$.
So we can write $A(s)=\exp(\x(s))^*B(s)$ with smooth paths
$\x:(-s_0,s_0)\to D^{1,p}(\d)\cap(\ker\rd_{A_0})^\perp$
and $B:(-s_0,s_0)\to L_{A_0}=\cL\cap S_{A_0}(\ep)$ such that
$\x(0)=0$ and $B(0)=A(0)$.

Next, by Proposition~\ref{prp lag slice}
we have a trivialization of $\rT L_{A_0}$ near $A_0$,
\[
\P:\; \begin{aligned}
L_{A_0} \times \rT_{A_0}L_{A_0}
\; & \to & S_{A_0} \\
(B,\b) \;& \mapsto & \P(B)\b ,
\end{aligned}
\]
such that $\P(B):\rT_{A_0}L_{A_0}\to\rT_{B}L_{A_0}$ is an isomorphism
for all $B$ sufficiently $L^2$-close to $A_0$.
We use this to write
$\pd_s B(s)=\P(B(s))\b(s)$ with a smooth path $\b : (-s_0,s_0)\to\rT_{A_0}L_{A_0}$.
Now we have
$$
\pd_s A(s) 
= \exp(-\x(s))\bigl(\P(B(s))\b(s)\bigr)\exp(\x(s)) 
 + \rd_{A(s)}\bigl( \exp(-\x(s))\pd_s \exp(\x(s)) \bigr)
$$
and hence by $\x(0)=0$ and $B(0)=A(0)$
\begin{align*}
\pd_s A(0) 
&= \P(A(0))\b(0)  + \rd_{A(0)}\pd_s\x(0) ,\\
\pd_s^2 A(0) 
&= 
[\pd_s B(0)+\pd_sA (0),\pd_s\x(0)] 
+ \rT_{A(0)}\P (\pd_s B(0)) \b(0)
+ \P(A(0))\pd_s\b(0)  \\
&\quad 
+ \rd_{A(0)}\bigl( \pd_s^2\x(0) - \pd_s\x(0)\pd_s\x(0) \bigr) .
\end{align*}
Note here that the last two terms in $\pd_s^2 A(0)$ lie in $\rT_{A(0)}\cL$,
as does $\pd_s A(0)$. 
So the symplectic form on $\pd_s A(0)$ and $\pd_s^2 A(0)$
simplifies as follows. (From now on all calculations will be at $s=0$.)
\begin{align*}
&\int_\S \la \pd_s A \wedge \pd_s\pd_s A \ra  \\
&= 
\int_\S \biggl( \la [ \pd_s A \wedge (\pd_s B+\pd_s A)],\pd_s\x \ra
+ \la \pd_s A \wedge \rT_{A}\P (\pd_s B) \b \ra \\
&\leq 
 C_0 \|\pd_s A\|_{L^2} \|\pd_s B + \pd_s A \|_{L^2}
\|\pd_s\x\|_{W^{1,2}} 
+ \|\pd_s A\|_{L^2} \|\rT_{A}\P\| \|\pd_s B\| \|\b\| \\
&\leq  C_{\rT\cL} \|\pd_s A\|_{L^2(\S)}^3 .
\end{align*}
Here we used Corollary~\ref{cor div curl}~(ii) and the fact that
$\rd_{A}\pd_s A=\pd_s F_A=0$ as well as
$\rd_{A(0)}\pd_s B(0)=\pd_s F_B(0)=0$.
From Remark~\ref{rmk slice} we have 
$\|\rd_{A_0}\pd_s\x\|_{L^2}\leq C_1\|\pd_s A\|_{L^2}$,
and so since $\pd_s\x\in(\ker\rd_{A_0})^\perp$
$$
\|\pd_s \x\|_{W^{1,2}} 
\leq C_2 \|\rd_{A_0}\pd_s\x\|_{L^2} 
\leq C_1 C_2 \|\pd_s A\|_{L^2} .
$$
As a consequence we obtain for 
$\pd_s B(0)= \pd_s A(0) - \rd_{A(0)}\pd_s\x(0)$
$$
\|\pd_s B\|_{L^2} 
\leq  \| \pd_s A\|_{L^2}  
+ \|\rd_{A_0}\pd_s\x\|_{L^2} 
+ \|[(A-A_0),\pd_s\x]\|_{L^2} 
\leq C_3 \|\pd_s A\|_{L^2} .
$$
Since $B$ is a path in the finite dimensional manifold $L_{A_0}$,
all norms on $\pd_s B\in\rT_{B}L_{A_0}$ are equivalent.
%
%
The same applies to the path $\b$ in $\rT_{A_0}L_{A_0}$.
So we dropped the subscripts from these norms and just note that
$\|\b\|\leq C_4\|\pd_s A\|_{L^2}$ since $\P(B)\b=\pd_s B$
and $\P(B)$ is an isomorphism that is uniformly invertible
for $B$ in a neighbourhood of $A_0$.
Finally, we used a uniform bound on 
$\rT_B\P: \rT_B L_{A_0} \times \rT_{A_0}L_{A_0} \to L^p(\S,\rT^*\S\otimes\cg)$
for $B$ in a neighbourhood of $A_0$.
\QED

Secondly, restating Lemma~\ref{lem ext i},
we show that $\cL\subset\cA^{0,p}(\S)$ is uniformly locally quasiconvex in the $L^2$-metric, despite possibly not being a topological submanifold.

\begin{lem}  \label{lem ext}
There are universal constants $C_\cL$ and $\d_\cL>0$ such that for all 
$A_1, A_2\in\cL$ with $\|A_1-A_2\|_{L^2}\leq\d_\cL$ there exists a path 
$\tA:[0,1]\to \cL$ with $\tA(0)=A_1$, $\tA(1)=A_2$, and
\begin{equation} \label{eqn est}
 \|\pd_s \tA (s) \|_{L^2} \leq C_\cL \|A_1-A_2\|_{L^2}  \qquad\forall s\in[0,1].
\end{equation}
\end{lem}

Before embarking on the proof we should remark that this lemma would be
a mere corollary of Theorems~\ref{thm short}, \ref{thm local slice}, and
Proposition~\ref{prp lag slice}
if we knew that the local slice map is continuous in the sense that
it provides $u_1^*A_1,u_2^*A_2\in S_{A_0}$ in the local slice of a
fixed nearby $A_0$ such that 
$\|u_1^*A_1-u_2^*A_2\|_{L^2}\leq C_0 \|A_1-A_2\|_{L^2}$.
In the present proof we replace this unknown continuity
by our knowledge of the bundle structure of $\cL$.\\

\noindent
{\bf Proof of Lemma~\ref{lem ext}: }
As in \cite[Section~3]{W bubb} we choose standard generators
$\a_1,\dots,\a_{2g}:[0,1]\to\S$ of $\pi_1(\S,z)$ that coincide 
near $z$.
These extend to embeddings $\ta_i:[-1,1]\times[0,1]\to\S$ such that
the loops $\ta_i(\t,\cdot)$ are based at the same family 
$z:[-1,1]\to\S$.
This provides a family of holonomy morphisms $\r_{z(\t)} : \cL \to \rG^{2g}$ 
given by parallel transport along the loops $\ta_i(\t,\cdot)$.
Now \cite[Lemma~3.1]{W bubb} says that for all 
$A_1,A_2\in\cL$ there exists $\t\in[-1,1]$ such that
$$
{\rm dist}_{\rG^{2g}} \bigl( \r_{z(\t)}(A_1) \,,\,  \r_{z(\t)}(A_2) \bigr)
\;\leq\; C_1 \| A_1 - A_2 \|_{L^1(\S)} .
$$
By \cite[Lemma~4.3]{W Cauchy} the images $\r_{z(\t)}(\cL)\subset\rG^{2g}$
are smooth submanifolds. 
Moreover, $\cL$ is a $\cG^{1,p}_{z(\t)}(\S)$-bundle over the compact 
quotients $\cL/\cG^{1,p}_{z(\t)}(\S)\cong\r_{z(\t)}(\cL)=:M_\t\subset\rG^{2g}$.
We fix a finite cover $M_0=\bigcup_{j=1}^N \cU_0^j$ 
such that there exist smooth local sections $\p_j:B_1\to\cL$ over 
the closed unit ball $B_1\subset\R^m$ inducing diffeomorphisms
$\r_{z(0)}\comp\p_j:B_1\to\cU_0^j\subset\rG^{2g}$.
Since the $\p_j$ are smooth over the compact $B_1$, the maps
$\r_{z(\t)}\comp\p_j:B_1\to M_\t\subset\rG^{2g}$ 
will also be local diffeomorphisms for small variations of $z(\t)$.
So there is a uniform constant $C_2$ such that 
for all $\t\in[-1,1]$ and $v_1,v_2\in B_1$
$$
\|v_1-v_2\|_{\R^m} \leq 
C_2 \, {\rm dist}_{\rG^{2g}} 
\bigl( \r_{z(\t)}(\p_j(v_1)) \,,\, \r_{z(\t)}(\p_j(v_2)) \bigr).
$$
Moreover, we can assume that the $\cU_\t^j:=\r_{z(\t)}(\p_j(B_1))\subset\rG^{2g}$ 
provide a cover of $M_\t$ for all $\t\in[-1,1]$.

Now given $A_1, A_2\in\cL$ with $\|A_1-A_2\|_{L^2}\leq\d$ sufficiently small
we find $\t\in[-1,1]$ such that
$\r_{z(\t)}(A_1),\r_{z(\t)}(A_2)\in\rG^{2g}$ are so close that they
lie in the same chart $\cU_\t^j$ for some $j$.
Then the $v_i:=(\r_{z(\t)}\comp\p_j)^{-1}\r_{z(\t)}(A_i)\in B_1$
satisfy
$$
\|v_1-v_2\|_{\R^m} \leq C_1 C_2 \| A_1 - A_2 \|_{L^1(\S)} .
$$
Moreover, we will have $A_i=u_i^*(\p_j(v_i))$ for some $u_i\in\cG^{1,p}_{z(\t)}(\S)$.
From this we can construct a first part of the required path.
$\tA(s):=u_2^*(\p_j(v(s))$ with $v(s):=(1-s) v_1 + s v_2)$
connects $\tA(1)=A_2$ to a connection $\tA(0)=u^*A_1$ 
which is gauge equivalent to $A_1$ by $u=u_1^{-1} u_2\in\cG^{1,p}(\S)$.
The length of this path is bounded by
$$
\|\pd_s \tA (s) \|_{L^2} 
= \|(\rT_{v(s)}\p_j) (v_2-v_1) \|_{L^2} 
\leq C_3 \|v_1-v_2\|_{\R^m} 
\leq C_1 C_2 C_3 \|A_1-A_2\|_{L^1} .
$$
So it remains to prove the lemma for gauge equivalent connections
$A_1=A,u^*A\in\cL$ with
\begin{align*}
\|A-u^*A\|_{L^2} &\leq \|A_1 - A_2\|_{L^2} + \| A_2 - u^* A_1 \|_{L^2}    \\
&\leq \|A_1 - A_2\|_{L^2} +  C_1C_2C_3 \| A_1 - A_2 \|_{L^1} 
\leq C_4 \| A_1 - A_2 \|_{L^2} \leq C_4 \d_\cL
\end{align*}
Note that the claim of the lemma for these is preserved under gauge transformation of $A$.
Now we need to digress for a moment to produce an open cover of $\cL/\cG^{1,p}(\S)$ by neighbourhoods on which we can achieve uniform constants.
Given any base point $A_0\in\cL$, Proposition~\ref{prp lag slice} provides $\ep_{A_0}>0$
such that $L_{A_0}:=\cL\cap S_{A_0}(\ep_{A_0})$ is a smooth submanifold and hence locally connected (though with respect to the $L^p$-topology).
So we can find $\delta_{A_0}>0$ such that any $A\in \cL\cap S_{A_0}$ with
$\|A-A_0\|_{L^p}\leq \delta_{A_0}$ lies in the same connected component of $L_{A_0}$ 
as $A_0$.
Since $\cL/\cG^{1,p}(\S)$ is compact, we can now assume without loss of 
generality that $A\in S_{A_0}(\ep)$ lies in the local slice of one
of finitely many base points $A_0\in\cL$, and moreover $\|A-A_0\|_{L^p}\leq \delta_{A_0}$.

Next, we can pick the constants $\ep_{A_0}>0$ and $\d_\cL>0$ such that
$\ep_{A_0}+C_4\d_\cL\leq\d$ for the $\d>0$ from Theorem~\ref{thm short}.
Then that theorem provides a path $v(s)^*A\in\cL$
from $v(0)^*A=A$ to $v(1)^*A$ such that 
$\rd_{A_0}^*(v(1)^{-1\;*}u^*A - A_0)=0$,
and whose length is bounded by
$$
 \bigl\| \pd_s \bigl( v(s)^*A \bigr) \bigr\|_{L^2}
\leq C_4 \bigl( 1+ \|A-A_0\|_{L^p} \bigr) \|u^*A-A\|_{L^2} 
\leq 2 C_5 \|u^*A-A\|_{L^2}.
$$
It remains to connect the endpoint $v(1)^*A$ to $u^*A$.
Equivalently we can connect $A$ to $w^*A$ (with $w=u\cdot v(1)^{-1}\in{\rm Stab}(A_0)$) 
for $A,w^*A\in S_{A_0}$ and
$$
\|A-w^*A\|_{L^2} 
= \|v(1)^*A-u^*A\|_{L^2}
\leq (1+2C_5) \|A - u^*A\|_{L^2} 
\leq C_6 \| A_1 - A_2 \|_{L^2}
\leq C_6\d_\cL .
$$ 
Since $w\in{\rm Stab}(A_0)$ we have $\|w^*A - A_0\|_{L^p}=\|A - A_0\|_{L^p}\leq\delta_{A_0}$, so by construction both $A$ and $w^*A$ lie in the connected component of $A_0$ in the
finite dimensional manifold $L_{A_0}=\cL\cap S_{A_0}(\ep_{A_0})$. 
They can thus be connected by a geodesic in $L_{A_0}\subset\cL$ 
whose length (and speed) is bounded linearly by $\|A-w^*A\|_{L^2}$.

If we first reparametrize the three separate paths above so that 
they are constant near the ends (and their slope in the interior 
is at most doubled), then the concatenated path is smooth and
satisfies the claimed bound on the derivative for all times.
\QED

\section{The Chern-Simons functional}
\label{sec CS}

Throughout this section we consider a gauge invariant Lagrangian submanifold 
$\cL\subset\cA^{0,p}$ as in Definition~\ref{dfn L}.
We moreover assume that the quotient space 
$\cL/\cG^{1,p}_z(\Sigma)$ is connected and simply connected for some (and hence every) 
base point set $z\subset\Sigma$.
The aim of this section is to define a local Chern-Simons functional 
for short arcs with endpoints on $\cL$ and establish an isoperimetric inequality.

We consider a smooth path \hbox{$A:[0,\pi]\to\cA(\S)$} with $L^2$-close 
endpoints $A(0),A(\pi)\in\cL$.
Then Lemma~\ref{lem ext} provides a continuous and piecewise smooth path 
${\tA:S^1\to\cA(\S)}$ with $\tA|_{[0,\pi]}\equiv A$ and $\tA([\pi,2\pi])\subset\cL$
such that 
\begin{equation}\label{ext cond}
\|\pd_\p \tA (\p) \|_{L^2} \leq C_\cL \|A(0)-A(1)\|_{L^2} 
\qquad\forall\p\in[\pi,2\pi] .
\end{equation}
We pick any such path to define the local Chern-Simons functional for $A$ 
by the usual Chern-Simons functional on $S^1\times\S$ for the extended connection $\tA$,
\begin{align}
\CS(A) &:= - \half\int_0^{2\pi} \int_\S \la \tA \wedge \pd_\p \tA \ra \;\dph .
\label{CS int}
\end{align}
A different choice of the extension path $\tA:[\pi,2\pi]\to\cL$ would change $\CS(A)$
by the value of the Chern-Simons functional on a loop $B:S^1\to\cL$.
This value however is invariant under homotopies: Let $B:[0,1]\times S^1 \to \cL$ be
continuous and piecewise smooth, then $\pd_\p B, \pd_s B\in\rT_B\cL$ for almost all
$(s,\p)$ and thus
\begin{align*}
\CS(B(1,\cdot))-\CS(B(0,\cdot))
&= \half\int_0^1 \frac\pd{\pd s} 
\int_{S^1} \int_\S \la B(s,\p) \wedge \pd_\p B(s,\p) \ra \;\dph\;\ds \\
&= - \int_0^1
\int_{S^1} \o\bigl(\pd_\p B(s,\p)\,,\,\pd_s B(s,\p)\bigr) \;\dph\;\ds 
\;=\; 0 .
\end{align*}
So the local Chern-Simons functional in (\ref{CS int}) is well defined up to
the additive subgroup $\CS(\pi_1(\cL))\subset\R$.
By our assumption on the base space $\cL/\cG^{1,p}_z(\S)$ of the bundle
(\ref{bundle}), the fundamental group $\pi_1(\cL)$ is generated by loops 
$u^*A_0$ for $A_0\in\cL$ and $u:S^1\to\cG(\S)$. 
For these we have $\CS(u^*A_0)=\CS(A_0)-4\pi^2\deg u \in 4\pi^2\Z$, and thus
$\CS(\pi_1(\cL))=4\pi^2\Z$.
With this we will see that our local Chern-Simons functional is in fact real valued
when restricted to sufficiently short paths, and it moreover satisfies an isoperimetric
inequality.

\begin{lem} {\bf (Isoperimetric inequality)} \label{lem iso} \\
There is $\ep>0$ such that for all smooth paths ${A:[0,\pi]\to\cA(\S)}$ with 
endpoints ${A(0),A(\pi)\in\cL}$ and $\int_0^\pi \|\pd_\p A\|_{L^2(\S)}\leq\ep$ 
the local Chern-Simons functional (\ref{CS int}) is well defined and
satisfies
$$
| \CS(A) | \;\leq\; \half (1+\pi C_\cL)^2 
\left( \int_0^\pi \bigl\| \pd_\p A \bigr\|_{L^2(\S)} \,\dph \right)^2 .
$$
\end{lem}

\noindent
{\bf Proof of Lemma~\ref{lem iso} :}\;
Let $A:[0,\pi]\to\cA(\S)$ be a smooth path with ${A(0),A(\pi)\in\cL}$
and $\int_0^\pi \|\pd_\p A\|_{L^2(\S)}\leq\ep$, where $\ep>0$ will be fixed
later on.
Consider any extending path $\tA:S^1\to\cA(\S)$ 
such that $\tA|_{[0,\pi]\times\S}\equiv A$, $\tA([\pi,2\pi])\subset\cL$, and
(\ref{ext cond}) holds.
With this we have

\begin{align*}
2\bigl|\CS(A)\bigr|
&= \biggl|\int_0^{2\pi}\int_\S \la  \tA \wedge \pd_\p \tA \ra \, \dph \biggr| \\
&= \biggl| \int_0^{2\pi}\int_\S \la \Bigl( \tA(0) + \int_0^\p \pd_\p \tA(\th)
 \,\rd\th \Bigr) \wedge \pd_\p \tA(\p) \ra \, \dph \biggr| \\
&= \biggl| \int_0^{2\pi}\int_0^\p \int_\S \la \pd_\p \tA(\th) 
           \wedge \pd_\p \tA(\p) \ra \, \rd\th \, \dph \biggl| \\
&\leq \left( \int_0^{2\pi} \bigl\| \pd_\p \tA(\p)  \bigr\|_{L^2(\S)} \, \dph  \right)^2 \\
&\leq  \left( \int_0^\pi \bigl\| \pd_\p A \bigr\|_{L^2(\S)} \,\dph 
            +  \pi C_\cL \| A(0)-A(\pi) \|_{L^2(\S)} \right)^2 \\
&\leq (1+\pi C_\cL)^2 \left( \int_0^\pi \bigl\| \pd_\p A \bigr\|_{L^2(\S)} \,\dph \right)^2 
\;\leq\; (1+\pi C_\cL)^2\ep^2.
\end{align*}
If we choose $\ep>0$ small enough, then this implies that our choice of extension path
will always yield values $\CS(A)\in [-\pi^2 , \pi^2]$.
If we change the path $\tA$, then as seen before
the Chern-Simons functional will change by a multiple of $4\pi^2$.
This cannot lead to another value in the interval $[-\pi^2 , \pi^2]$, hence the
value of $\CS(A)$ is uniquely determined by the condition (\ref{ext cond}) on the
extensions.
\QED

\subsubsection*{The Chern-Simons functional for handle body Lagrangians}

For a Lagrangian submanifold $\cL=\cL_H$ that arises from a handle body $H$
with $\pd H=\Sigma$ an alternative definition of the Chern-Simons functional
is given in \cite[Section~4]{W bubb}.
For $A:[0,\pi]\to\cA(\S)$ with $L^2$-close endpoints $A(0),A(\pi)\in\cL_H$
we defined
\begin{align*}
\CS_H(A) &:= 
- \half\int_0^{\pi} \int_\S \la A \wedge \pd_\p A \ra \;\dph 
- \tfrac 1{12} \int_H \winner{A_H(0)}{[A_H(0)\wedge A_H(0)]} \\
&\qquad\qquad\qquad\qquad\qquad\qquad
+ \tfrac 1{12} \int_H \winner{A_H(\pi)}{[A_H(\pi)\wedge A_H(\pi)]} ,
\end{align*}
where flat extensions $A_H\in\cA_{\rm flat}(H)$ are chosen such that 
$A_H|_\Sigma=A|_\Sigma$ and 
$\|A_H(0)-A_H(\pi)\|_{L^3(H)} \leq C_H \|A(0)-A(\pi)\|_{L^2(\Sigma)}$
with some uniform constant $C_H$.
Here we show that $\CS_H$ agrees with the more general $\CS$ given by 
(\ref{CS int}).

The extensions $A_H(0),A_H(\pi)\in\cA_{\rm flat}(H)$ determine
a path $\tA:[\pi,2\pi]\to\cL_H$ uniquely up to homotopy,
by requiring that $\tA$ extends to a path
$\tilde A_H:[\pi,2\pi]\to\cA_{\rm flat}(H)$ with the given endpoints
$\tA_H(\pi)=A_H(\pi)$ and $\tA_H(2\pi)=A_H(0)$.
This is since $\cA_{\rm flat}(H)$ is a bundle over the simply 
connected base $\SU(2)\times \dots\times \SU(2)$, 
whose fibre $\cG_z(H)$ is also simply connected 
(since $H$ retracts onto its $1$-skeleton and $\pi_2(\SU(2))=0$).
With this path we indeed obtain $\CS(A)=\CS_H(A)$ since
\begin{align*}
\tfrac 12 \int_\S \la \tA \wedge \pd_\p \tA \ra 
&=
\tfrac 12\int_H \Bigl( 
\la \rd_{\tA_H}\tA_H \wedge \pd_\p \tA_H \ra 
- \la \tA_H \wedge \rd_{\tA_H} \pd_\p \tA_H \ra 
\Bigr) \\
&=
\tfrac 1{12}\; \frac\rd{\rd\ph} \int_H \la [\tA_H\wedge\tA_H] \wedge \tA_H \ra  
\end{align*}
Here $F_{\tA_H}=0$, so
$\rd_{\tA_H}\tA_H=[\tA_H\wedge\tA_H]$ and
$\rd_{\tA_H} \pd_\p \tA_H=\pd_\ph F_{\tA_H}=0$.

\section{Bubbling for ASD instantons with general Lagrangian boundary conditions}
\label{sec bubb}

In this section, based on the results of Sections~\ref{sec Lag} and~\ref{sec CS},
we extend the compactness results of \cite{W bubb} for anti-self-dual connections with Lagrangian boundary conditions \eqref{bvpi}
to the more general class of Lagrangian boundary conditions $\cL$ 
introduced in Section~\ref{sec CS}.

\begin{thm}
Let $\cL\subset\cA^{0,p}$ be a gauge invariant Lagrangian submanifold 
as in Definition~\ref{dfn L} for some $p>2$, 
and suppose that the quotient space 
$\cL/\cG^{1,p}_z(\Sigma)$ is connected and simply connected for some base point set $z\subset\Sigma$.
Then the energy quantization Theorem \cite[Thm.1.2]{W bubb}
and the removable singularity Theorem \cite[Thm.1.5]{W bubb}
continue to hold with $\cL_Y$ replaced by $\cL$.
This finishes the proof of compactness for moduli spaces of \eqref{bvpi},
as already claimed in Theorem \cite[Thm.7.2]{SW}.
\end{thm}

Rather than copying statements and proofs we provide new
results and indicate how these can replace the (few but crucial) 
arguments in \cite{W bubb} that are based on the special form of the
Lagrangians.
Roughly, we only need to replace \cite[Lemma~1.6]{W bubb}
and the definition of the local Chern-Simons functional.

As in \cite{W bubb} we fix a metric of normal type $\ds^2+\dt^2+g_{s,t}$ 
on $D\times\S$, where $D:=B_{r_0}(0)\cap\H^2$ is the $2$-dimensional 
half ball of radius $r_0>0$ and centre~$0$.
We consider a connection $\X\in\cA(D\times\S)$ that solves the boundary value problem
\begin{equation}\label{bvp}
F_\X + * F_\X = 0, \qquad \qquad
\X|_{(s,0)\times\Si} \in\cL \quad\forall s\in[-r_0,r_0] .
\end{equation}
In the proof of the energy quantization
\cite[Theorem~1.2]{W bubb} we only need to replace
the estimate for the normal derivative in \cite[Lemma~2.3]{W bubb}
by the following result.
We express the connection in the splitting $\Xi=A+\Phi\ds+\Psi\dt$
and denote by $B_s=\pd_s A -\rd_A\P$ one of the curvature components.

\begin{lem} \label{lem norm}
There is a constant $C$ (varying continuously with the 
metric of normal type in the $\cC^2$-topology) 
such that for all solutions $\X\in\cA(D\times\Si)$ of (\ref{bvp})
\begin{align*}
-\tfrac\pd{\pd t}\bigr|_{t=0} \bigl\| F_\X \bigr\|_{L^2(\Si)}^2 
&\;\leq\; C \bigl( \bigl\| B_s \bigr\|_{L^2(\Si)}^2 
+ \bigl\| B_s \bigr\|_{L^2(\Si)}^3 \bigr) .
\end{align*}
\end{lem}
\Pr
As in \cite[Lemma 2.3]{W bubb} we have
\begin{align*}
- \tfrac 14 \tfrac\pd{\pd t}\bigr|_{t=0} \bigl\| F_\X \bigr\|_{L^2(\Si)}^2 
&\leq \Bigl( C \bigl\| B_s \bigr\|_{L^2(\Si)}^2 
            - \int_\Si \la \nabla_s B_s \wedge B_s \ra \Bigr)\Bigr|_{t=0} .
\end{align*}
The estimate for this normal derivative can be checked in any gauge at a fixed 
$(s_0,0)\in D\cap\pd\H^2$.
We choose a gauge with $\P\equiv 0$ and hence $B_s=\pd_s A$.
Then $\X|_{(\cdot,0)\times\S}=A(\cdot,0)$ is a path in $\cL$
to which Lemma~\ref{lem lin ext} applies.
So we calculate
\begin{align*}
 - \int_\Si \la \nabla_s B_s \wedge B_s \ra  
\;=\;  \int_\S \la \pd_s A \wedge \pd_s\pd_s A \ra  
\;\leq\; C_{\rT\cL} \bigl\| \pd_s A \bigr\|_{L^2(\Si)}^3 
\;=\; C_{\rT\cL} \bigl\| B_s \|_{L^2(\Si)}^3 .
\end{align*}
The constant $C_{\rT\cL}$ does not depend on the metric and the
constant $C$ above varies as in \cite{W bubb}.
\QED

For the removal of singularity Theorem \cite[Thm.1.5]{W bubb}
we denote punctured half balls by
$$
D^*_r:= B_r(0)\setminus\{0\} \,\cap\, \H^2 , \qquad\qquad
D^* := D^*_{r_0}
$$
and consider a finite energy solution $\X\in\cA(D^*\times\S)$ of (\ref{bvp}).
Using polar coordinates $r\in(0,r_0]$, $\p\in[0,\pi]$ on $D^*$ we can always
choose a gauge $\X=A+R\dr$ with no $\dph$-component.
Then the energy function is
\begin{align}\label{energy b}
\E(\r) &:= \half \int_{D^*_\r\times\S} |F_\X|^2  
&= \int_0^\r \int_0^\pi \bigl( \|F_A\|_{L^2(\S)}^2 
+ r^{-2}\|\pd_\p A\|_{L^2(\S)}^2 \bigr) \,r\,\dph\,\dr .
\end{align}
In the proof of \cite[Theorem~1.5]{W bubb} we only need to 
replace the isoperimetric inequality \cite[Lemma~4.1~(ii)]{W bubb}.
This is reproven in the following lemma.

\begin{lem} \label{lem isoper}
Let $\X\in\cA(D^*\times\S)$ satisfy (\ref{bvp}) and $\E(r_0)<\infty$,
and suppose that it is in the gauge $\X=A+R\dr$ with $\P\equiv 0$.
Then there exists $0<r_1\leq r_0$ such that for all $r\leq r_1$
$$
\E(r) \;=\; - \CS(A(r,\cdot)) 
\;\leq\; \half (1+\pi C_\cL)^2
\left( \int_0^\pi \bigl\| \pd_\p A(r,\p)  \bigr\|_{L^2(\S)} \, \dph  \right)^2 
\;\leq\; \half\b^{-1} r \,\dot\E(r)
$$
and hence $\E(r)\leq C r^{2\b}$ with some constants $C$ and $\b>0$.
\end{lem}

Its proof uses \cite[Lemma~5.4]{W bubb}, restated below, 
which directly generalizes to the present case 
(using Lemma~\ref{lem norm} in the proof instead
of the normal estimate from \cite[Lemma~2.3]{W bubb}).

\begin{lem} \label{lem mean curv}
There exist constants $C$ and $\ep>0$ such that the following holds.
Let $\X$ be a smooth connection on $D^*\times\S$ that satisfies (\ref{bvp}).
Suppose that $\E(2r)\leq\ep$ for some $0<r\leq \half r_0$, then for all 
$\p\in(0,\pi)$
$$
\|F_\X(r,\p)\|_{L^2(\S)}^2 \leq C r^{-2} \E(2r) 
\qquad\text{and}\qquad
\|F_\X(r,\p)\|_{L^\infty(\S)}^2 \leq C (r\sin\p)^{-4} \E(2r) .
$$
\end{lem}

\noindent
{\bf Proof of Lemma~\ref{lem isoper}:}\;
Let $0<\r\leq r_1$ (where $0<r_1<r_0$ will be fixed later on), 
then by assumption $\E(\r)\leq\E(r_0)$ is finite, i.e.\ it exists as the limit
\begin{align*}
\E(\r) &= \lim_{\d\to 0} \half \int_{( D_\r\setminus D_\d ) \times \S} |F_\X|^2 .
\end{align*}
We aim to express this as the difference of a functional at $r=\r$ and at $r=\d$. 
The straightforward approach picks up additional boundary terms
on $\{\p=0\}$ and $\{\p=\pi\}$. We eliminate these by gluing in $[\pi,2\pi]\times\S$
and extending the connection by paths in $\cL$.

More precisely, $(D_\r\setminus D_\d)\times\S$ is diffeomorphic (with preserved orientation)
to $[\d,\r]\times[0,\pi]\times\S$, and we glue in $[\d,\r]\times [\pi,2\pi]\times\S$ 
to obtain the smooth $4$-manifold
$[\d,\r]\times S^1\times\S$ which has the boundary components
$S^1 \times \S$ at $r=\r$ and $S^1 \times \bar\S$ at $r=\d$.
Next, $A(\cdot,0)$ and $A(\cdot,\pi)$ are smooth paths in $\cL$, and for sufficiently
small $r_1>0$ they will automatically lie in the same connected component 
(see (\ref{?!}) below).
So we can pick a smooth family of extension paths $\tA : [\d,\r]\times S^1 \to \cA(\S)$
with $\tA|_{[\d,\r]\times[0,\pi]}\equiv A$ and 
$\tA([\d,\r],[\pi,2\pi])\subset\cL$.
We also extend the function $R:[\d,\r]\times[0,\pi]\times\S\to\cg$ 
to a smooth function $\tR:[\d,\r]\times S^1\times\S\to\cg$.

These extensions match up to a $W^{1,\infty}$-connection 
$\tX:=\tA + \tR\dr$ on $[\d,\r]\times S^1\times\S$.
Note that the extension over $[\d,\r]\times [\pi,2\pi]\times\S$ will not contribute 
to the energy expression for the instanton, that is
\begin{align}
&\half \int_\d^\r \int_\pi^{2\pi} \int_\S \la F_\tX \wedge F_\tX \ra  \nonumber\\
&=
\int_\d^\r \int_\pi^{2\pi} \int_\S \Bigl( \la F_\tA \,,\, \pd_\p \tR \ra
           + \la (\pd_r\tA -\rd_\tA \tR ) \wedge \pd_\p \tA \ra \Bigr) \;\dph\;\dr
\;=\;0 . \label{zero}
\end{align}
This uses the fact that $\tA\in\cL$, hence
$F_\tA\equiv 0$ and 
$\int_\S \la (\pd_r\tA -\rd_\tA R)\wedge \pd_\p A\ra = 0$ since this is the 
symplectic form on tangent vectors to the Lagrangian~$\cL$.

We will choose the extension paths $\tA(\cdot,[\pi,2\pi])$
such that for all ${\d\leq r\leq\r}$ the functional $\cC(\tA(r,\cdot))$
given by the right hand side of (\ref{CS int}) with this extension path 
equals to the local Chern-Simons functional $\CS(A(r,\cdot))$.
For this purpose let $\bep>0$ be the constant from Lemma~\ref{lem mean curv} and choose 
$0<r_1\leq \half r_0$ such that $\E(2r_1)\leq\bep$. Then for all ${0<r\leq r_1}$
\begin{align} \label{?!}
\biggl(\int_0^\pi \bigl\| \pd_\p A(r,\p)  \bigr\|_{L^2(\S)} \, \dph \biggr)^2
&\leq \pi \int_0^\pi \bigl\| \pd_\p A(r,\p)  \bigr\|_{L^2(\S)}^2 \, \dph \\
&\leq \tfrac \pi 2 \int_0^\pi r^2 \bigl\| F_\X(r,\p) \bigr\|_{L^2(\S)}^2 \,\dph  
\;\leq\; C \E(2r) . \nonumber
\end{align}
Now choose $r_1>0$ even smaller such that $C\E(2r_1)\leq\min(\pi^2,\ep^2)$
with $\ep>0$ from Lemma~\ref{lem iso}, and such that $A(r,0)$ and $A(r,\pi)$ 
automatically lie in the same connected component of $\cL$ for $r\leq r_1$.
Then the lemma applies to $A(r,\cdot)$ for all $0<r\leq r_1$.
In particular, since $\r\leq r_1$, we can choose the family of extension paths
to end at $\tA(\r,\cdot)$ such that 
$$
\cC(\tA(\r,\cdot)):= 
- \half\int_0^{2\pi} \int_\S \la \tA(\r,\p) \wedge \pd_\p \tA(\r,\p) \ra \;\dph
\;=\;\CS(A(\r,\cdot)) .
$$
Moreover we know that for all $r\in [\d,\r]$ the path $A(r,\cdot)$ is sufficiently
short for the local Chern-Simons functional $\CS(A(r,\cdot))$ to be defined
and take values in $[-\pi^2,\pi^2]$.
Now $\cC(\tA(r,\cdot))$ is a smooth function of $r\in[\d,\r]$
whose values might differ from $\CS(A(r,\cdot))$ by multiples of $4\pi^2$.
We have equality at $r=\r$ and hence by continuity for all $r\in[\d,\r]$ as claimed.
Thus we actually express the energy of $\tX$ in terms
of the local Chern-Simons functional 
\begin{align} 
\half \int_{( D_\r\setminus D_\d ) \times \S} |F_\X|^2 
&= - \half \int_{[\d,\r]\times S^1\times\S} \la F_\tX \wedge F_\tX \ra \nonumber\\
&= - \cC(\tA(\r,\cdot)) + \cC(\tA(\d,\cdot))  \label{energy CS}  \\
&= - \CS(A(\r,\cdot)) + \CS(A(\d,\cdot)) .\nonumber
\end{align}
Here we have $F_\tX\wedge F_\tX = -|F_\X|^2 \dvol$ on $(D_\r\setminus D_\d ) \times \S$
and $\int \la F_\tX\wedge F_\tX \ra = 0$ 
on $[\d,\r]\times [\pi,2\pi]\times\S$ by (\ref{zero}).
Now by Lemma~\ref{lem iso} 
\begin{align*}
\bigl| \CS(A(r,\cdot)) \bigr|
&\leq \half (1+\pi C_\cL)^2 
\left( \int_0^\pi \bigl\| \pd_\p A(r,\p) \bigr\|_{L^2(\S)} \,\dph \right)^2 \\
&\leq \frac{\pi(1+\pi C_\cL)^2}4 
                \int_0^\pi r^2 \bigl\| F_\X(r,\p) \bigr\|_{L^2(\S)}^2 \,\dph  .
\end{align*}
As $r\to 0$ this expression converges to zero by Lemma~\ref{lem mean curv}.
Thus for all ${0<\r\leq r_1}$ 
$$
\E(\r) 
\;=\; - \CS(A(\r,\cdot))
\;\leq\; \tfrac 14 \b^{-1} \int_0^\pi \r^2 \bigl\| F_\X(\r,\p) \bigr\|_{L^2(\S)}^2 \,\dph  
\;\;=\; \half\b^{-1} \r \, \dot \E(\r) 
$$
with $\b=\pi^{-1}(1+\pi C_\cL)^{-2}>0$.
By integration this implies $\E(r) \leq C r^{2\b} $ for all $0<r\leq r_1$.
\QED

\bibliographystyle{plain}

\end{document}